\documentclass[11pt]{article}
\usepackage{amsmath,amsxtra}
\usepackage{tabularx}
\usepackage{enumerate}
\usepackage{amssymb}
\usepackage{amscd}
\usepackage{amsthm}
\usepackage[latin1]{inputenc}
\usepackage{t1enc}
\makeatletter
\DeclareRobustCommand\sfrac[1]{\@ifnextchar/{\@sfrac{#1}}%
                                            {\@sfrac{#1}/}}
\def\@sfrac#1/#2{\leavevmode\raise.5ex
         \hbox{$\m@th\mbox{\fontsize\sf@size\z@
                           \selectfont#1}$}\kern-.1em
         /\kern-.1em\lower.25ex
          \hbox{$\m@th\mbox{\fontsize\sf@size\z@
                            \selectfont#2}$}}
\makeatother
\def\1e{\frac{1}{\e}}

\def\L{{\mathcal L}}
\def\half{1/2}
\def\e{\varepsilon}
\def\/{\, | \,}
\def\n{\nabla}

\def\ee{\epsilon}
\def\div{\delta}
\def\dr{{\mathcal D}}
\def\d{\text{ d}}

\def\ox{\otimes}
\def\K{{\mathcal K}}

\def\RE{{\mathbf R}} 
\def\NA{{\mathbb N}} 

\def\car{{\mathbf 1}}
\def\I{{\mathcal I}}

\def\S{{\mathcal S}}

\def\car{{\mathbf 1}}
\def\egaldef{\stackrel{def}{=}}
\def\egalnot{\stackrel{not}{=}}

\def\P{{\text P}}

\def\D{{\mathbb D}}

\def\<{\langle}
\def\>{\rangle}
\def\({\Bigl(}
\def\){\Bigr)}
\def\pth{p_{t,t+h}}
\def\ks{\K^*_1}

\def\xth{\frac{Z_t+Z_{t+h}}{2}}
\def\D{{\mathbb D}}
\newcommand{\esp}[1]{{\text E}\left[{#1}\right]}
\newcommand{\II}[3]{\text{\tiny{(\textsc{#1})-}}\kern-3pt\int #3(s)\, #2 dB_H(s)} 

\newcommand{\trace}{\operatorname{trace}}
\newcommand{\Dom}{\operatorname{Dom}}
\newcommand{\Id}{\operatorname{Id}}

\newcommand{\Hol}{\operatorname{Hol}}

\def\DD{{\mathcal D}}
\def\DDp{{\mathcal D}^\prime}
{\theoremstyle{plain}
\newtheorem{prop}{Proposition}[section]
\newtheorem{thm}{Theorem}[section]
\newtheorem{cor}{Corollary}[section]
\newtheorem{lemma}{Lemma}[section]
\newtheorem{defn}{Definition}[section]
}
{\theoremstyle{definition}
\newtheorem{hyp}{Hypothesis}
\newtheorem{notation}{Notation}}

{\theoremstyle{remark}
\newtheorem{rem}{Remark}[section]
\newtheorem{example}{\textbf{Example}}
}
\makeatletter
\newenvironment{example*}[1][\examplename]{\par
  \normalfont
  \topsep6\p@\@plus6\p@ \trivlist
  \item[\hskip\labelsep\itshape
    \textbf{Example #1}\@addpunct{.}]\ignorespaces
}{%
  $\sqsupset$\endtrivlist
}

\newcommand{\examplename}{\textbf{Example}}
\makeatother
\title{Stochastic Integration with respect to Volterra processes}
\author{L. Decreusefond} \date{}
\begin{document}
\maketitle
\begin{abstract}
  We construct the basis of a stochastic calculus for so-called
  Volterra processes, i.e., processes which are defined as the
  stochastic integral of a time-dependent kernel with respect to a
  standard Brownian motion. For these processes which are natural
  generalization of fractional Brownian motion, we construct a
  stochastic integral and show some of
  its main properties: regularity with respect to time and kernel,
  transformation under an absolutely continuous change of probability,
  possible approximation schemes and Itô formula.
\end{abstract}  
\section{Introduction}
\label{sec:introduction}
In the past few years, more than twenty papers have been devoted to the
definition of a stochastic integral with respect to fractional
Brownian motion or other ``related'' processes, see for instance
\cite{decreusefond99_4} and references therein.
Remind that  fractional Brownian process of Hurst index $H\in
(0,1),$ denoted by $B^H,$ is the unique centered
Gaussian process whose covariance kernel is given by 
 \begin{equation*}
    R_H(s,t)=\esp{B_s^HB_t^H}\egaldef
    \frac{V_H}{2}\(s^{2H}+t^{2H}-|t-s|^{2H}\)
  \end{equation*}
  where
\begin{equation*}
  V_H\egaldef \frac{\Gamma(2-2H)\cos (\pi H)}{\pi H (1-2H)}.
\end{equation*}
Among other properties, this process has $1/H$-finite variation and a finite generalized  covariation of
order 4 for $H>1/4,$ (see \cite{russo02} for the definition),  has H\"older continuous
trajectories of any order less than $H$ and has the following
representation property:
\begin{equation}
\label{eq:23}
  B^H(t)=\int_0^t K_H(t,s) \d B_s,
\end{equation}
where $B$ is a one dimensional standard Brownian motion and $K$ is
deterministic kernel with an intricate expression (see
\cite{decreusefond96_2}). Therefore, a ``related'' process means altogether a process
with finite $p$-variation, called a process with rough paths in
\cite{coutin01,lyons98}, or  a process with H\"older
continuous sample-paths as in \cite{feyel96,zaehle99} and also a
process of the form (\ref{eq:23}) with a general kernel as in
\cite{alos00,coutin02,decreusefond02}.

This is the last track that we will follow here. Our present work,
which is the expanded version of \cite{decreusefond02}, differs from
the other two papers \cite{alos00,coutin02} in two ways. First, the
method to define the stochastic integral is different. In these two
papers, the kernel is regularized, if needed, to obtain a
semi-martingale. The second step is then to  use the
classical theory of stochastic integration and then pass to the limit
after a stochastic integration by parts in the sense of the Malliavin
Calculus. We here use an approach based on convergence of discrete
sums. It should be already noted  that for smooth integrands, their
 notion of integral and ours coincide. The other difference is to be
 found in the kind of hypothesis put on $K.$ In \cite{alos00,coutin02},
 hypothesis are made on the regularity of the function $K(t,s)$
 itself. We
 here work with assumptions on the linear map $f\mapsto \int
 K(t,s)f(s)\d s.$ Properties of $K(t,s)$ and $Kf$ are, of course,
 intimately related but we think that working with the latter gives
 more insight on the underlying problems.
 
 In Section \ref{sec:prelim}, we recall basic definitions and
 properties of deterministic fractional calculus. In Section
 \ref{sec:malliavin}, we introduce the class of processes, named
 Volterrra processes, that we will study. We then give a few
 properties of their sample-paths. In Section \ref{sec:approx}, we deal
 with  a Stratonovitch-like definition of the stochastic integral
 with respect to Volterra processes. Section \ref{sec:regul} is
 devoted to the time regularity of the previously constructed integral
 and in Section \ref{sec:ito}, we establish an It\^o formula. In the
 last section, we show how the Stratonovith integral is related to a
 Skorohod-like integral and how a It\^o-like process constructed from such an
 integral is modified through an absolutely continuous change of probability.
 
\section{Preliminaries}
\label{sec:prelim}
This section is only devoted to the presentation of the tools of
deterministic fractional calculus we shall use in the sequel.
For $f\in \L^1([0,1];\ dt),$ (denoted by $\L^1$ for short) the left
and right fractional integrals of $f$ are defined by~:
        \begin{align*}
          (I_{0^+}^{\gamma}f)(x) & \egaldef
          \frac{1}{\Gamma(\gamma)}\int_0^xf(t)(x-t)^{\gamma-1}dt\ ,\ 
          x\ge
          0,\\
          (I_{1^-}^{\gamma}f)(x) & \egaldef
          \frac{1}{\Gamma(\gamma)}\int_x^1f(t)(t-x)^{\gamma-1}dt\ ,\ 
          x\le 1,
        \end{align*}
        where $\gamma>0$ and $I^0_{0^+}=I^0_{1^-}=\Id.$ For any
        $\gamma\ge 0$, any $f\in \L^p$ and $g\in \L^q$ where
        $p^{-1}+q^{-1}\le \gamma$, we have~:
\begin{equation}
  \label{int_parties_frac}
  \int_0^1 f(s)(I_{0^+}^\gamma g)(s)\ ds = \int_0^1 (I_{1^-}^\gamma 
f)(s)g(s)\ ds.
\end{equation}
The Besov-Liouville space $I^\gamma_{0^+}(\L^p)\egalnot \I_{\gamma,p}^+$ is
usually equipped with the norm~:
\begin{equation}
\label{normedansIap}
  \|  I^{\gamma}_{0^+}f \| _{ \I_{\gamma,p}^+}=\| f\|_{\L^p}.
\end{equation}
Analogously, the Besov-Liouville space $I^\gamma_{1^-}(\L^p)\egalnot \I_{\gamma,p}^-$ is
usually equipped with the norm~:
\begin{equation*}
  \| I^{-\gamma}_{1^-}f \| _{ \I_{\gamma,p}^-}=\|  f\|_{\L^p}.
\end{equation*}
We then have the following continuity results (see
\cite{feyel96,samko93})~:
\begin{prop}
\label{prop:proprietes_int_RL}
  \begin{enumerate}[i.]
  \item \label{inclusionLpLq} If $0<\gamma <1,$ $1< p <1/\gamma,$ then
    $I^\gamma_{0^+}$ is a bounded operator from $\L^p$ into $\L^q$
    with $q=p(1-\gamma p)^{-1}.$
  \item  For any $0< \gamma <1$ and any $p\ge
    1,$ $\I_{\gamma,p}^+$ is continuously embedded in $\Hol(\gamma-
    1/p)$ provided that $\gamma-1/p>0.$ $\Hol(\nu)$ denotes the space
    of H\"older-continuous functions, null at time $0,$ equipped with
    the usual norm.
  \item  For any $0< \gamma< \beta <1,$ $\Hol
    (\beta)$ is compactly embedded in $\I_{\gamma,\infty}.$
\item \label{semigroupe} By $I^{-\gamma}_{0^+},$ respectively
  $I^{-\gamma}_{1^-},$ we mean the inverse map of $I^{\gamma}_{0^+},$
  respectively $I^{\gamma}_{1^-}.$ The relation
  $I^{\gamma}_{0^+}I^{\beta}_{0^+}f=I^{\gamma+\beta}_{0^+}f$ holds
  whenever $\beta >0,\ \gamma+\beta >0$ and $f\in \L^1.$
\item For $\gamma p>1,$ the spaces $\I_{\gamma,p}^+$ and
  $\I_{\gamma,p}^{-}$ are canonically isomorphic. We will thus use the
  notation $\I_{\gamma,p}$ to denote any of this spaces. This property
  isn't any more
  true for $\gamma p>1,$ see Lemma  \ref{lem:inclusion_restriction}
  and text below Definition \ref{def:strat-integr}. 
  \end{enumerate}
\end{prop}
We now define the Besov-Liouville spaces of negative order and
 show that they are in duality with Besov-Liouville of positive order
 (it is likely that this exists elsewhere in the literature but we
 have not found any reference so far).
Denote by $\DD_+$ the space of ${\mathcal C}^\infty$ functions defined
on $[0,1]$ and such that $\phi^{(k)}(0)=0.$ Analogously, set  $\DD_-$ the space of ${\mathcal C}^\infty$ functions defined
on $[0,1]$ and such that $\phi^{(k)}(1)=0.$ They are both equipped
with  the  projective topology induced
by the semi-norms $p_k(\phi)=\sum_{j\le k}\lVert \phi^{(j)} \rVert_\infty.$ Let
$\DD^\prime_+,$ resp. $\DD^\prime_-,$  be  their strong topological dual. It is straightforward
that $\DD_+$ is stable by $I^\gamma_{0^+}$ and $\DD_-$ is stable
$I^\gamma_{1^-},$ for any $\gamma\in \RE.$ Hence, guided by (\ref{int_parties_frac}),
we can define the fractional integral of any distribution (i.e., an
element of $\DDp_-$ or $\DDp_+$):
\begin{align*}
\text{ For }T\in \DDp_-;\   I^\gamma_{0^+}T: \ \phi \in\DD_- & \mapsto <T,\,
  I^\gamma_{1^-}\phi>_{\DDp_-,\DD_-},\\
\text{ For }T\in \DDp_+;\     I^\gamma_{1^-}T: \ \phi \in\DD_+ & \mapsto <T,\,
  I^\gamma_{0^+}\phi>_{\DDp_+,\DD_+}.
\end{align*}
We introduce now our Besov spaces of negative order by 
\begin{defn}
  For $\gamma > 0$ and $r>1,$ $\I_{-\gamma,r}^+$
  (resp. $\I_{-\gamma,r}^-$) is the space of 
  distributions such that $I^\gamma_{0^+}T$  (resp. $I^\gamma_{1^-}T$
  ) belongs to $\L^r.$ The
  norm of an element $T$ in this space is the norm of
  $I^\gamma_{0^+}T$ in $\L^r$ (resp. of  $I^\gamma_{1^-}T$).
\end{defn}
\begin{thm}
\label{thm:cardual}
For $\gamma > 0$ and $r>1,$ the dual space of $\I_{\gamma,r}^+$
(resp. $\I_{\gamma,r}^-$)
is canonically isometrically isomorphic to
$I^{-\gamma}_{1^-}(\L^{r^*})$ (resp. $I^{-\gamma}_{0^+}(\L^{r^*})$,) 
where $r^*=r(r-1)^{-1}.$
\end{thm}
\begin{proof}
 Let $T$ be in $\DDp_+,$ we have:
 \begin{align*}
   \sup_{\phi:\lVert
     \phi \rVert_{\I_{\gamma,r}^+}=1}|<T,\, \phi>| & = \sup_{\psi:\lVert
     \psi \rVert_{\L^r}=1}|<T,\,I_{0^+}^\gamma \phi>|\\
& =\sup_{\psi:\lVert
     \psi \rVert_{\L^r}=1}|<I_{1^-}^{\gamma}T,\, \phi>|
 \end{align*}
hence   by the Hahn-Banach theorem,
  \begin{equation*}
   T \in (\I_{\gamma,r}^+)^\prime  \Longleftrightarrow \sup_{\phi:\lVert
     \phi \rVert_{\I_{\gamma,r}^+}=1}|<T,\, \phi>| < \infty
 \Longleftrightarrow I_{1^-}^{\gamma} T\in \L^{r^*},
  \end{equation*}
and $\lVert T \rVert_{(\I_{\gamma,r}^+)^\prime} =\lVert T
\rVert_{I^{-\gamma}_{1^-}(\L^{r^*})}.$ The same reasoning also holds
for $(\I_{\gamma,r}^-)^\prime.$
\end{proof}
\begin{thm}
\label{thm:integrationdssoboneg}
  For $\beta\ge \gamma \ge 0$ and $r>1,$  $I_{1^-}^\beta$ is continuous from
  $\I_{-\gamma,r}^-$ into $\I_{\beta-\gamma,r}^-.$ 
\end{thm}
\begin{proof}
Since $T$ belongs to $\I_{-\gamma,r}^-=(\I_{\gamma,r^*})^\prime,$ we
have:
\begin{equation*}
|  <I_{1^-}^\beta T,\, \phi>| =|  < T,\,I_{0^+}^\beta \phi>| 
 \le c
\lVert I_{0^+}^\beta \phi\rVert_{\I_{\gamma,r^*}}
=c \lVert I_{0^+}^{\beta-\gamma} \phi\rVert_{\L^{r^*}}.
\end{equation*}
Thus, $I_{1^-}^\beta T$ is a continuous linear form on
$\I_{\gamma-\beta,r^*}^+$ and thus belongs to the dual of this space
which, according to the previous theorem, is exactly $\I_{\beta-\gamma,r}^-.$ 
\end{proof}
For $\eta >0$ and $p\in [1,+\infty),$ the Slobodetzki space $\S_{\eta,p}$ is the closure of ${\mathcal
  C}^1$ functions with respect to the semi-norm:
\begin{equation*}
  \label{eq:10}
  \lVert f\rVert_{\S_{\eta,p}}^p=\iint_{[0,1]^2}
  \frac{|f(x)-f(y)|^p}{|x-y|^{1+p\eta}}\, dx\, dy,
\end{equation*}
For $\eta=0,$ we simply have $\S_{0,p}=L^p([0,1]).$
We then have the following continuity results (see
\cite{feyel96,zaehle99})~:
\begin{prop}
  \begin{enumerate}[i.]
    \item  For any $0< \gamma <1$ and any $p\ge
    1,$ $\S_{\gamma,p}$ is continuously embedded in $\Hol(\gamma-
    1/p)$ provided that $\gamma-1/p>0.$ $\Hol(\nu)$ denotes the space
    of H\"older-continuous functions, null at time $0,$ equipped with
    the usual norm. 

For $0<\gamma<1/p,$ $\S_{\gamma,p}$ is compactly embedded in
$L^{p(1-\gamma p)^{-1}}([0,1]).$ Moreover, if $p=2,$ the embedding of
$\S_{\gamma,p}$ into $L^2([0,1])$ is Hilbert-Schmidt.
\item \label{besov_embed}
It is proved in \cite{feyel96} that for $1\ge a>b>c>0$ that we the
following embeddings are continuous (even compact)
\begin{equation}
  \label{eq:17}
 \S_{a,p}\subset \I_{b,p}^+\subset \S_{c,p}.
\end{equation}
  \item  For any $0< \gamma< \beta <1,$ $\Hol
    (\beta)$ is compactly embedded in $\S_{\gamma,\infty}.$
\item \label{besov_nesting} Let $a>0,$ $1< p\le q<\infty.$ Suppose
  $b=a-1/p+1/q>0.$ Then $\S_{a,p}$ is continuously embedded in
  $\S_{b,q},$ see \cite{adams75}.
    \end{enumerate}
\end{prop}

One of the key property we shall use, is this result due to Tambaca \cite{tambaca01}.
\begin{lemma}
\label{lem:tambaca}
Let $r,s\in [0,\half)$ and let $t=r+s-\half\ge 0.$ For $f\in \S_{s,2},$
$g\in \S_{r,2},$ the product $fg$ belongs to $\S_{t,2}$ and we have:
  \begin{equation*}
    \lVert fg \rVert_{\S_{t,2}}\le c  \lVert f \rVert_{\S_{r,2}} \lVert g \rVert_{\S_{s,2}}.
  \end{equation*}
\end{lemma}
From this Lemma and the embeddings of Eqn. (\ref{eq:17}), we have: 
\begin{cor}
\label{cor:tambaca}
  Let $r,s\in [-\infty,\half)$ and let $t<r+s-\half.$ For $f\in \I_{s,2},$
$g\in \I_{r,2},$ the product $fg$ belongs to $\I_{t,2}$ and we have:
  \begin{equation*}
    \lVert fg \rVert_{\I_{t,2}}\le c  \lVert f \rVert_{\I_{r,2}} \lVert g \rVert_{\I_{s,2}}.
  \end{equation*}
\end{cor}
We will need a similar result in the simpler situation where $r$ is
greater than $1/2.$
\begin{lemma}
  Let $r>1/2,$ for $f$ and $g$ in $\I_{r,2},$ we have
  \begin{equation}
    \label{eq:tambaca_facile}
    \lVert fg\rVert_{\S_{r,2}}\le c \lVert f\rVert_{\S_{r,2}}\lVert g\rVert_{\S_{r,2}}.
  \end{equation}
\end{lemma}
\begin{proof}
  Since $r>1/2,$ $f$ and $g$ are continuous and $\lVert
  f\rVert_\infty\le c \lVert f\rVert_{\S_{r,2}}.$ The same holds for
  $g.$
Thus,
\begin{align*}
  \lVert fg\rVert_{\S_{r,2}}^2&\le \iint_{[0,1]^2}
  \Bigl(\frac{|f(x)|^2(g(x)-g(y))^2}{|x-y|^{1+2r}}+\frac{|g(y)|^2(f(x)-f(y))^2}{|x-y|^{1+2r}}\Bigr)\d x\d y\\
& \le c \Bigl(\lVert
  f\rVert_\infty^2\lVert g\rVert_{\S_{r,2}}^2+\lVert
  g\rVert_\infty^2\lVert f\rVert_{\S_{r,2}}^2\Bigr),
\end{align*}
and the result follows.
\end{proof}

One could probably work with only one family of spaces (i.e., either
$\I_{\alpha,p}$ or $\S_{\alpha,p}$) but depending on the properties,
some are easier to verify in the setting of Riemann-Liouville spaces  and some
 in the setting of Slobodetzki spaces, see for instance the property
 below.
 \begin{lemma}
\label{lem:inclusion_restriction}
   Let $\gamma >\tilde{\gamma}>1/2$ and $f \in \S_{\gamma,2}$ then
   $(f-f(t))\car_{[0,t]}$ belongs to $\S_{\tilde{\gamma},2}.$ 
 \end{lemma}
 \begin{proof}
   First note that $f$ is $(\gamma-1/2)$-H\"older continuous thus that $f-f(t)$ is well
   defined. Moreover,
   \begin{multline*}
\iint_{[0,1]^2}
  \frac{|(f(x)-f(t))\car_{[0,t]}-(f(y)-f(t))\car_{[0,t]}|^2}{|x-y|^{1+2\tilde{\gamma}}}\, dx\, dy \\
 =   \iint_{[0,t]^2}
  \frac{|f(x)-f(y)|^2}{|x-y|^{1+2\tilde{\gamma}}}\, dx\, dy + 2 \iint_{[0,t]\times[0,1]}
  \frac{|f(x)-f(t)|^2}{|x-y|^{1+2\tilde{\gamma}}}\, dx\, dy\\
\le \lVert f\rVert_{\tilde{\gamma},2}^2(1+2\iint_{[0,t]\times[t,1]}
  \frac{|x-t|^{2\gamma-1}}{|x-y|^{1+2\tilde{\gamma}}}\, dx\, dy)
\le c \lVert f\rVert_{\tilde{\gamma},2}^2.
   \end{multline*}
 \end{proof}
%
%
\section{Volterra processes}
\label{sec:malliavin}
Consider that we are given a deterministic Hilbert-Schmidt linear map,
$K,$ satisfying:
\begin{hyp}
  \label{A1}
  There exists $\alpha >0$ such that $K$ is continuous, one-to-one,
  from $\L^2([0,1])$ into $\I_{\alpha+1/2,2}.$ Moreover, $K$ is
  triangular, i.e., for any $\lambda \in [0,1],$ the set ${\mathcal
    N}_\lambda =\{ f: \ f(t)=0 \text{ for } t\le \lambda\}$ is
  invariant by $K.$
\end{hyp}
\begin{rem}
  Since $K$ is Hilbert-Schmidt from $\L^2([0,1])$ into itself, there exists a measurable kernel
  $K(.,.)$ such that
  \begin{equation*}
    Kf(t)=\int_0^1 K(t,s)f(s)\d s.
  \end{equation*}
  The triangularity of $K$ is equivalent to $K(t,s)=0$ for $s>t,$
  i.e.,
  \begin{equation*}
    Kf(t)=\int_0^t K(t,s)f(s)\d s.
  \end{equation*}
\end{rem}
Consider now the kernel $R(t,s)$ defined by
\begin{displaymath}
R(t,s):= \int_0^{t\wedge s} K(t,r)K(s,r)\d r.
\end{displaymath}
The map associated to $R,$ i.e., $Rf(t)=\int_0^1 R(t,s)f(s)\d s,$ is
equal to $KK^*$ and for any $\beta_1,\ldots,\beta_n$ any $t_1,\ldots,
t_n,$ we have
\begin{displaymath}
  \sum_{i,j} \beta_i\beta_j R(t_i,t_j)=\int K^*(\sum \beta_j
  \epsilon_{t_j})(s)^2\d s\ge 0,
\end{displaymath}
so that $R(t,s)$ is a positive kernel and we can speak of the centered
Gaussian process of covariance kernel $R.$ Let $X$ be this process and
be the subject of our study.
\begin{lemma}
  The process $X$ has a modification with a.s. continuous
  sample-paths.
\end{lemma}
\begin{proof}
  We have
  \begin{align*}
    \esp{(X_t-X_s)^2}&=\int_0^t K(t,r)^2 \d r+\int_0^s K(s,r)^2 \d
    r-2\int_0^{t\wedge s} K(t,r)K(s,r) \d r\\
&= K(K(t,.)-K(s,.))(t)-K(K(t,.)-K(s,.))(s)\\
&\le c |t-s|^\alpha \Bigl(\int_0^1 (K(t,r)-K(s,r))^2\d r
\Bigr)^{1/2}.
  \end{align*}
Expanding the square in the last integral, we get the right hand side
of the first equation, thus 
\begin{equation*}
   \esp{(X_t-X_s)^2}^{1/2}\le c |t-s|^\alpha.
\end{equation*}
Kolmogorov Lemma entails that $X$ has a modification with H\"older
continuous sample paths of any order less than $\alpha.$
\end{proof}
We thus now work on the Wiener space $\Omega={\mathcal
  C}_0([0,1];{\mathbf R}),$ the Cameron-Martin space is
$H=K(\L^2([0,1]))$ and $P,$ the probability on $\Omega$ under which the
canonical process, denoted by $X,$ is a centered Gaussian process of
covariance kernel $R.$ The norm of $h=K(g)$ in $H$ is the norm of $g$
in $\L^2([0,1]).$

A mapping $\phi$ from $\Omega$ into some separable Hilbert space $X$
is called cylindrical if it is of the form
$\phi(w)=\sum_{i=1}^df_i(\<v_{i,1},w\>,\cdots,\<v_{i,n},w\>) x_i$
where for each $i,$ $f_i\in {\mathcal C}_0^\infty (\RE^n, \RE)$ and
$(v_{i,j},\, j=1\ldots n)$ is a sequence of $\Omega^*$ such that
$(\tilde{v}_{i,j},\, j=1\ldots n)$ (where $\tilde{v}_{i,j}$ is the
image of $v_{i,j}$ under the injection $\Omega^\star\hookrightarrow
\L^2([0,1])$ ) is an orthonormal system of $\L^2([0,1]).$ For such a
function we define $\n \phi$ as
$$
\n\phi(w)=\sum_{i,j=1} \partial_j
f_i(\<v_{i,1},w\>,\cdots,\<v_{i,n},w\>){\tilde{v}}_{i,j}\otimes x_i.
$$
From the quasi-invariance of the Wiener measure
\cite{ustunel_book}, it follows that $\n$ is a closable operator on
$L^p(\Omega;X)$, $p\geq 1$, and we will denote its closure with the
same notation. The powers of $\n$ are defined by iterating this
procedure. For $p>1$, $k\in \NA$, we denote by $\D_{p,k}(X)$ the
completion of $X$-valued cylindrical functions under the following
norm
$$
\|\phi\|_{p,k}=\sum_{i=0}^k \|\n^i\phi\|_{L^p(\Omega;X\otimes
  \L^2([0,1])^{\otimes i})}\,.
$$
\begin{rem}
Note that the Sobolev spaces $\S_{\alpha,p}$ enjoy the useful property
of $p$-admissibility (after \cite{feyel91}) and thus for any $0< \gamma <1$ and any $p\ge
    1,$ the spaces $\D_{p,k}(\S_{\alpha,p})$ and
    $\S_{\alpha,p}(\D_{p,k})$ are isomorphic.  
\end{rem}

The divergence, denoted $\div$ is the adjoint of $\n$: $v$ belongs
to $\Dom_p \div$ whenever for any cylindrical $\phi,$
\begin{equation*}
  |\esp{\int_0^1 u_s \n_s\phi\d s}|\le \, c \lVert \phi\rVert_{L^p}
\end{equation*}
and for such a process $v,$ $$\esp{\int_0^1 u_s \n_s\phi\d s}=\esp{
  \phi\, \div u}.$$
It is easy to show (see \cite{decreusefond96_2})
that $\{B_t:=\div(\car_{[0,t]}),\, t\ge 0\}$ is a standard Brownian
motion such that $\div u=\int u_s\d B_s$ for any square integrable
adapted processes $u$ and which satisfies
$$X_t=\int_0^t K(t,s) \d B_s.$$
Moreover, $B$ and $X$ have the same
filtration. In view of the last identity and because $K$ is lower
triangular, we decided to name such a process, a Gaussian Volterra
process. The analysis of processes of the same kind where $B$ is
replaced by a jump processes is the subject of our current
investigations with N. Savy.
\begin{example}
  \label{fbm_levy}
The first example is the so-called L\'evy fractional Brownian
motion of Hurst index $H$, defined as 
\begin{displaymath}
 \frac{1}{\Gamma(H+1/2)} \int_0^t (t-s)^{H-1/2}\d B_s.
\end{displaymath}
This amounts to say that $K=I_{0^+}^{H+1/2},$ thus that hypothesis
\ref{A1} and \ref{A2} are immediately satisfied, with $\alpha=H,$ in view of  the semi-group properties of
 fractional integration.
\end{example}
\begin{example}
  \label{fbm}
The other classical example is the fractional Brownian motion with
stationary increments of Hurst index $H,$ for which 
   \begin{equation}
  \label{defdekh}
K(t,s)=  K_H(t,r):=\frac{(t-r)^{H- \frac{1}{2}}}{\Gamma(H+\frac{1}{2})}
F(\frac{1}{2}-H,H-\frac{1}{2}, H+\frac{1}{2},1- 
\frac{t}{r})1_{[0,t)}(r).
\end{equation}
The Gauss hyper-geometric function $F(\alpha,\beta,\gamma,z)$  (see
\cite{nikiforov88}) is the analytic continuation on ${\mathbb
  C}\times {\mathbb C}\times {\mathbb C} \backslash \{-1,-2,\ldots
\}\times \{z\in {\mathbb C}, Arg |1-z| < \pi\}$ of the power series
 \begin{displaymath}
    \sum_{k=0}^{+ \infty} \frac{(\alpha)_k(\beta)_k}{(\gamma)_k 
k!}z^k,
 \end{displaymath}
and 
\begin{displaymath}
  (a)_0=1 \text{ and } (a)_k \egaldef 
\frac{\Gamma(a+k)}{\Gamma(a)}=a(a+
1)\dots (a+k-1).
\end{displaymath}
We know from \cite{samko93} that 
$K_H$ is an isomorphism from $\L^2 ([0,1])$ onto 
                $\I_{H+1/2,2}^+$ and 
                \begin{align*}
                        K_Hf &= 
                I_{0^+}^{2H}x^{\half-H}I_{0^+}^{\half-H}x^{H-\half}f \ \text{ 
                        for } H \le  \half,\\
                        K_Hf &= 
                I_{0^+}^1x^{H-\half}I_{0^+}^{H-\half}x^{\half-H}f \ \text{ 
                        for } H \ge \half.
                \end{align*}
It follows easily that Hypothesis \ref{A1} and \ref{A2} are satisfied
with $\alpha=H.$
\end{example}
\begin{example}\label{ex:gfbm}
  Beyond these two well known  cases, we can investigate the case of $K(t,s)=K_{H(t)}(t,s)$
  for a deterministic function $H.$ This is the process studied in
  \cite{benassi99}. It seems interesting to analyze since statistical
  investigations via wavelets have shown that the local H\"older
  exponent of some real signals (in telecommunications) is varying
  with time and this situation can't be reflected with a model based
  on fBm since its H\"older regularity is everywhere equal to its
  Hurst index.
  \begin{lemma}
    For $f\in \L^2,$ for $H_1>H_2\ge \gamma>0,$ we have
    \begin{align}
\label{eq:ineqKH2}
|K_{H_2}f(s)-K_{H_2}f(t)|&\le c |t-s|^\gamma \lVert f\rVert_{\L^2},\\
\label{eq:ineqKH1}      |K_{H_1}f(s)-K_{H_2}f(s)|& \le c |H_1-H_2|
\lVert f\rVert_{\L^2},
    \end{align}
    where $c$ is a constant independent of $H_1,\, H_2$ and $f.$
  \end{lemma}
  \begin{proof}
Since $H_2$ is greater than $\gamma,$ $K_{H_2}f$ belongs to
$\I_{\gamma+1/2,2},$ and (\ref{eq:ineqKH2}) follows directly from the embedding of
    $\I_{\gamma+1/2,2}$ into $\Hol(\gamma).$ 

Another expression of the hypergeometric function is given by:
\begin{displaymath}
  F(a,b,c,z)=\frac{\Gamma(c)}{\Gamma(b)\Gamma(c-b)}\int_0^1
  t^{c-1}(1-t)^{c-b-1}(1-zt)^{-a}\d t.
\end{displaymath}
Classical and tedious computations show that for $H\in[h_1+\e,h_2-\e],$
\begin{displaymath}
  |\frac{d}{dH}K_H(t,s)|\le c_\e \sup_{H\in (H_1,H_2)}|K_H(t,s)|,
\end{displaymath}
where $c_\e=\sup_{t\in[0,1]}|t^\e\ln t|.$ It thus
entails that 
$$|K_{H_2}(t,s)-K_{H_1}(t,s)| \le c_\e \sup_{H\in (H_1,H_2)}|K_H(t,s)|
|H_2-H_1|.$$ Cauchy-Schwarz inequality yields to (\ref{eq:ineqKH1}).
  \end{proof}
  \begin{thm}
    Let $H$ belong to $\S_{1/2+\alpha,2}$ and be such that $\inf_t H(t)>1/2,$
    then $K(t,s)=K_{H(t)}(t,s)$ satisfies \ref{A1} for any  $\alpha<\inf_t H(t)-1/2.$
\end{thm}    
\begin{proof}
Let $f$ belong to $\L^2,$  set $\gamma=\inf_t H(t)$ and let $\alpha<\gamma-1/2.$
According to the previous lemma, we have
\begin{multline*}
  \lVert Kf\rVert_{\S_{1/2+\alpha,2}}^2  =\iint_{[0,1]^2}
  \frac{|K_{H(t)}f(t)-K_{H(s)}f(s)|^2}{|t-s|^{2+2\alpha}} \d t  \d s\\
 \shoveleft{ \le 2\iint_{[0,1]^2}
  \frac{|K_{H(t)}f(t)-K_{H(t)}f(s)|^2}{|t-s|^{2+2\alpha}} \d t \d s}\\
\shoveright{+  2\iint_{[0,1]^2}
  \frac{|K_{H(t)}f(s)-K_{H(s)}f(s)|^2}{|t-s|^{2+2\alpha}} \d t  \d s}\\
  \le c \,\lVert f\rVert_{\L^2}^2\iint_{[0,1]^2}
  \frac{|t-s|^{2\gamma}}{|t-s|^{1+2\alpha}}
  \d t \d s
+ c\,\lVert f\rVert_{\L^2}^2\iint_{[0,1]^2}
  \frac{|H(t)-H(s)|^2}{|t-s|^{2+2\alpha}} \d t \d s.
\end{multline*}
 The right-hand-side is finite by hypothesis and thus $K$ is
 continuous from $\L^2$ into $ \S_{1/2+\alpha,2}.$
\end{proof}
\end{example}
\section{Stratonovitch integral}
\label{sec:approx}
Starting
from scratch and trying  to define a stochastic integral with
respect to $X$ by a limit of a sequence of finite sums, we have two
main choices: Either we discretize $X$ (or more probably
$dX$) or we discretize $B$ (likely $dB$) and then derive a
discretization of $dX.$ The first approach yields two possibilities:
for a partition $\pi$ whose points are denoted by
$0=t_0<t_1<\ldots<t_n=T,$ we can consider
\begin{align} 
  \text{RS}_{\pi}(u)& =\sum_{t_i\in\pi} u(t_i) \Delta X_i \text{ or }\\
  \text{SS}_{\pi}(u)&=\sum_{t_i\in\pi}
  \frac{1}{\delta_i}(\int_{t_i}^{t_{i+1}}\!\!u(s)\d s)\  \Delta X_i,
\end{align}
where $\delta_i=t_{i+1}-t_i$ and $\Delta X_i=X(t_{i+1})-X(t_i).$ They
are both reminiscences of respectively Riemann and
Skorohod-Stratonovitch sums as defined in \cite{nualart.book}. 

In the other approach, we first
linearize $B$ and then look at the approximation of $X$ it yields
to. Let
\begin{equation*}
  B^\pi(t)=B(t_i)+\frac{1}{\delta_i}\Delta B_i (t-t_i) \text{ for } t\in[t_i,t_{i+1}),
\end{equation*}
and
\begin{align*}
  X^\pi(t)&
  =\sum_{t_i\in\pi}\frac{1}{\delta_i}\int_{t_i}^{t_{i+1}}K(t,s)\, ds \ 
  \Delta B_i \\
  &=\sum_{t_i\in\pi}\frac{1}{\delta_i} K(\car_{[t_i,t_{i+1}]})(t)
  \Delta B_i .
\end{align*}
It follows that it is  reasonable to consider
\begin{equation*}
R_T^\pi(u):=\sum_{t_i\in\pi}\frac{1}{\delta_i} \left\{\int_0^T u(t)
  \frac{d}{dt}K(\car_{[t_i,t_{i+1}]})(t) \d t\right\}\ \Delta B_i ,
\end{equation*}
under the additional hypothesis that
  for any $b>0,$ the function $K(\car_{[0,b]})$ is differentiable with a
  square integrable derivative.
For $u$ sufficiently smooth in the sense of the calculus of
variations, we have
\begin{multline*}
  R_T^\pi(u)= \div\Bigl(\sum_{t_i\in\pi}\frac{1}{\delta_i} \int_0^T u(t)
  \frac{d}{dt}K(\car_{[t_i,t_{i+1}]})(t) \, dt
  \car_{[t_i,t_{i+1}]}\Bigr)\\
+ \sum_{t_i\in\pi}\frac{1}{\delta_i} \int_{t_i}^{t_{i+1}}
\int_0^T \n_ru(t)\frac{d}{dt}K(\car_{[t_i,t_{i+1}]})(t) \d t\d r.
\end{multline*}
Using $\K^*_T,$ the formal adjoint of $\K:=I^{-1}_{0^+}\circ K$ on
$\L^2([0,T]),$ we have 
\begin{multline}
\label{eq:def_de_RpiT}   R_T^\pi(u)= \div\Bigl(\sum_{t_i\in\pi}\frac{1}{\delta_i}
   \int_{t_i}^{t_{i+1}} \kern-3pt \K^*_Tu(t) \d t\Bigr)
+ \sum_{t_i\in\pi}\frac{1}{\delta_i} \iint\limits_{[t_i,t_{i+1}]^2} \K^*_T(\n_ru)(t) \d t\d r.
\end{multline}
We now recognize the Skorohod-Stratonovitch sum associated to the
standard Brownian motion $B$ and to the integrand $\K^*_Tu.$ For the
sequel to be meaningful, we need to assume that the map $\K$
exists. This is guaranteed for $\alpha\ge 1/2,$ since
$\I^+_{\alpha+1/2,2}$ is embedded in the set of absolutely continuous
functions with square integrable derivative, but for $\alpha <1/2,$ we
need to introduce an additional hypothesis.
  \begin{hyp}
  \label{A2}
  We assume that for any $T\in [0,1],$ the map $\K=I_{0^+}^{-1}\circ K$ is a densely defined,
  closable operator from $\L^2([0,T])$ into itself and that its domain
  contains a dense subset, $\dr,$ stable by the maps $p_T,$ for any $T\in
  [0,1],$ where $p_Tf\equiv f\car_{[0,T)}.$  We denote by $\K^*_T$ its
  adjoint in $\L^2([0,T]).$ We assume furthermore that $\K^*_1$ is continuous from
  $\I_{1/2-\alpha,p}^{1^-}$ into $\L^2([0,T]),$ for any $p\ge 2.$
\end{hyp}
\begin{rem}
  In the preceding examples, $\dr$ may be taken to $\I_{(1/2-\alpha)^+,2}.$
\end{rem}
\begin{rem}
  For the sake of simplicity, we will speak of the domains of $\K$ and
  $\K^*_T$ independently of the position of $\alpha$ with respect to
  $\half.$ It must be plain that for $\alpha>\half,$ $\Dom \K
  =\L^2([0,1])$ and $\Dom
  \K^*_T=\L^2([0,T]).$
\end{rem}
\begin{rem}
  Since $I^1_{1^-}(\varepsilon_t)=\car_{[0,t]},$ we have 
  \begin{equation*}
    \K^*(\car_{[0,t]})=K^*(\varepsilon_t)=K(t,.).
  \end{equation*}
This means that $\K^*_t$ is identical to the operator denoted by
$I^{K_H}_t$ in \cite{coutin02}.
\end{rem}
\begin{notation}
  For any $p\ge 1,$ we denote by $p^*$ the conjugate of $p.$ For any
  linear map $A,$ we denote by $A^*_{T},$ its adjoint in
  $\L^2([0,T]).$ We denote by $c$ any irrelevant constant appearing in
  the computations, $c$ may vary from one line to another.
\end{notation} 
\begin{defn}
\label{def:strat-integr}
Assume that Hypothesis \ref{A1} holds for $\alpha \ge 1/2.$ We say that $u$ is Stratonovitch integrable on $[0,T]$ whenever
the family $\text{R}^\pi_T(u),$ defined in (\ref{eq:def_de_RpiT}), converges in probability as $|\pi|$ goes to $0.$ In
this case the limit will be denoted by $\int_0^T u_s \circ \d X_s.$
\end{defn}
This definition could be theoretically extended to $\alpha <1/2$
but would be practically unusable. Indeed, as we shall see below, when $\alpha<1/2,$ the
convergence of
the second sum of $R_T^\pi(u)$ requires that $u$ belongs to 
$\I_{1+\eta-\alpha,2}$ for some $\eta>0$ and $\K^*_T$ to be continuous
from this space to a space of Holderian functions. Since
$1+\eta-\alpha-1/2>0,$ the two spaces
$I^{1+\eta-\alpha}_{0^+}(\L^2([0,T]))$ and
$I^{1+\eta-\alpha}_{T^-}(\L^2([0,T]))$ are not canonically isomorphic
(if $u$ belongs to the first one then $u(0)=0$ whereas when $u$
belongs to the latter, $u(T)=0$). 
We thus have to specify to which one $u$ belongs exactly. In view of the
example of the L\'evy fractional Brownian where
$\K^*_T=I^{H-1/2}_{T^-},$ it is more convenient to assume that $u$ belongs to
$I^{1+\eta-H}_{T^-}(\L^2([0,T]))$ and thus that $u(T)$ is equal to $0.$
That raises a problem because the restriction of an element of
$I^{1+\eta-H}_{T^-}(\L^2([0,T]))$ to a shorter interval, say $[0,S],$
does not belong  $I^{1+\eta-H}_{S^-}(\L^2([0,S]))$ so that, we can't see
$\int_0^S u(r)\circ \d X_r$ as $\int_0^T u(r)\car_{[0,S]}(r)\circ \d X_r.$

On the other hand, since $(u-u(S))\car_{[0,S]}$ belongs to
$I^{1+\eta-H}_{S^-}(\L^2([0,S]))$ as soon as $u$ belongs to
$I^{1+\eta-H}_{T^-}(\L^2([0,T])),$ it is reasonable to consider
$R^\pi_T(u-u(T)).$ For the limit to stay the same, we have to add the
term $u(T)X(T).$ Indeed, the well known relationship (see \cite{nualart.book,ustunel_book}) 
\begin{equation}
  \label{eq:deltaaxi}
  \div(a\xi)=a\div\xi - \int_0^1 \n_r a \xi(r) \d r,
\end{equation}
for $a\in \D_{2,1}$ and $\xi\in\L^2(\Omega\times[0,1]),$ entails that
\begin{equation}
  \label{eq:lien_rpit_rpi1pt}
  R^\pi_T(u)=R^\pi_T(u-u(T))+u(T)X^\pi(T).
\end{equation}
As a
conclusion, for $\alpha <1/2,$ the definitive definition is 
\begin{defn}[Definition for $\alpha <1/2$]
Assume that Hypothesis \ref{A1} and \ref{A2} hold for $\alpha <1/2.$ We say that $u$ is Stratonovitch integrable on $[0,T],$ whenever
the family $\text{R}_T^\pi(u-u(T))$ converges in probability as $|\pi|$ goes to $0.$ In
this case, we set
\begin{equation}
  \label{eq:9}
  \int_0^T u_s \circ \d X_s =\lim_{|\pi|\to 0}R^\pi_T(u-u(T)) +u(T)X(T).
\end{equation}
\end{defn}
In view of the preceding discussion, the following lemma will play a key role in the sequel.
\begin{lemma}
\label{lem:triangularite}
For $T\in (0,1],$ let $p_Tf$ denote the restriction of $f$ to $[0,T).$
For any $f\in \Dom \K^*_1,$ $f$ belongs to $\Dom\K^*_T,$  $p_T f$
belongs to $\Dom\K^*_1$ and we have
  \begin{equation}
   \label{eq:5}
p_T\K^*_1(p_Tf)\equiv \K^*_T(f) .
     \end{equation}
\end{lemma}
\begin{proof}
  Since $K$ is triangular, for $g\in\dr,$ $p_T g$
  belongs to $\Dom\K$ and $p_T Kg=p_TK(p_Tg)=Kp_T g.$ By derivation,
   it follows that $p_T\K g=p_T\K p_Tg=\K p_Tg,$ so
  that, for $f\in\Dom\K^*_1,$
\begin{align*}
  |\int_0^t f(s)\K g(s)\d s|& = |\int_0^1 (p_T f)(s)\K g(s)\d s  |\\
  & = |\int_0^1  f(s)(p_T \K g)(s)\d s  |\\
  &=|\int_0^1 f(s) \K(p_T g)(s)\d s| \\
  &\le c \lVert p_T g\rVert_{\L^2([0,1])}=c \lVert
  g\rVert_{\L^2([0,T])}.
\end{align*}
By density, this identity
remains true for $g\in \Dom\K,$ thus this means that $f$ belongs to $\Dom\K^*_T$ and that $p_T f$ belongs
to $\Dom \K^*_1.$

For $g\in \L^2([0,T])\cap\Dom \K,$ we denote by $\tilde{g}$ its
extension to $\L^2([0,T])$ defined by $\tilde{g}(s)=0$ whenever $s\ge
T.$ We have
\begin{align*}
  \int_0^T p_T\K^*_1p_T f(s)g(s)\d s& = \int_0^1 \K^*_1p_T f(s)
  p_T\tilde{g}(s)\d s\\
  & = \int_0^1 p_Tf(s) \K (p_T\tilde{g})(s)\d s\\
  &= \int_0^T f(s)\K g(s)\d s\\
  &=\int_0^T \K^*_Tf(s) g(s)\d s,
\end{align*}
where the last equality follows by the first part of the proof and the
definition of the adjoint of a linear map. Since $g$ can be arbitrary,
(\ref{eq:5}) follows by identification.
\end{proof}
\begin{thm}
\label{thm:existence_trace-}
Let $\alpha <1/2$ and $p\ge 2.$ Assume that Hypothesis \ref{A1} and \ref{A2} hold. Assume furthermore
that  there exists $\sigma>1/p$ and $\eta>0,$  such that  $\K^*_1$ is continuous from
$\I_{\sigma,p}^{1^-}$ into $\Hol(\eta).$ If $u$ belongs to
$\D_{p,1}(\I_{\sigma+\e,p}^{1^-}),$ for some $\e>0,$  then for any $T\in [0,1],$ there exists a measurable and
integrable process, denoted by $\tilde{D}_Tu$ such that, for any $s,$
any $0\le a < b<1,$
\begin{multline}
  \label{eq:13}
\esp{\int_a^b | \K^*_T(\n_r (u-u(T)))(s)-\tilde{D}_Tu(r)|^p\d r}
\\
\le
\,c\,\esp{\int_0^1  |s-r|^{p\eta}\lVert \n_r
u\rVert_{\I^{1^-}_{\sigma+\e,p}}^p\d r}.  
\end{multline}
Moreover,
  \begin{equation}
    \label{eq:6bis}
    \esp{\lVert \int_0^.\tilde{D}_Tu(r)\d r\rVert_{\I_{1,p}^+}^p}\le c \lVert
    u\rVert_{\D_{p,1}(\I^{1^-}_{\sigma+\e,p})}^p.
  \end{equation}
\end{thm}
\begin{proof}
Since $\sigma>1/p,$ $u$ is continuous and we can speak unambiguously of
$u(T).$  The assumed continuity of $\K^*_1$ entails that $\K^*_T(u-u(T))$ belongs to
  $\D_{p,1}(\Hol(\eta))$ and that 
    \begin{multline}
\label{eq:borne_nablakappastaru}
\esp{\int_a^b      |\n _r \K^*_T(u-u(T))(s)-\n _r
  \K^*_T(u-u(T))(\tau)|^p\d r}\\
\le c \,\esp{\int_0^1 
      |s-\tau|^{p\eta}\lVert \n_r
      u\rVert_{\I_{\sigma+\e,p}^{1^-}}^p\d r}.
    \end{multline}
Consider $(\rho_n,\, n\ge 1)$ a one-dimensional positive mollifier, we
can define $\P\otimes \d r$ a.s.,  $\tilde{D}_Tu(s)$ by 
\begin{equation*}
\label{eq:def_de_D_t}
  \tilde{D}_Tu(r)=\lim_{n\to\infty} \int_0^T
  \rho_n(\tau)\K^*_T(\n_{r}u)(\tau-r)\d \tau.
\end{equation*}
   Hence, $\tilde{D}_Tu(r)$ is measurable with respect to $(\omega,r)$
   and according to (\ref{eq:borne_nablakappastaru}), we have
   (\ref{eq:13}). Substituting $0$ to $s$ (\ref{eq:borne_nablakappastaru}), we get 
\begin{equation*}
  \esp{\int_0^T |\tilde{D}_Tu(r)|^p\d s}\le c \lVert
    u\rVert_{\D_{p,1}(\I_{\sigma,p}^{1^-})}^p.
\end{equation*}
This means that $\int_0^.\tilde{D}_Tu(s)\d s$
belongs to $\I_{1,p}^+$  and that (\ref{eq:6bis}) holds.
\end{proof}
\begin{example*}[\protect\ref{fbm_levy} cont'd]
  In this case, $\K^*_1=I_{1^-}^{H-1/2}$ is continuous from
  $\I_{\sigma,p}^{1^-}$ into $\I_{\sigma+\alpha-1/2,p}^+.$ This latter
  space is embedded in a space of H\"olderian functions provided that $\sigma>1/2-\alpha+1/p.$
\end{example*}
\begin{example*}[\protect\ref{fbm} cont'd]
  According to \cite{samko93},
  $\K^*_1=x^{1/2-H}I^{H-1/2}_{1^-}x^{H-1/2}$ and since since
  $2(1+H-1/2)=2H+1>1,$ we infer from \cite[Lemma 10.1]{samko93} 
  that $\K^*_1$  is continuous from  $\I_{\sigma,p}^{1^-}$ into $\I_{\sigma+\alpha-1/2,p}^+,$
 for any $\sigma\ge 0.$
\end{example*}
\begin{thm}
\label{thm:strato_int_inf}
Let $\alpha <1/2$ and $p\ge 2.$ Assume that Hypothesis \ref{A1} and \ref{A2} hold. Assume furthermore
that  there exists $\sigma>1/p$ and $\eta>0,$  such that  $\K^*_1$ is continuous from
$\I_{\sigma,p}^{1^-}$ into $\Hol(\eta).$ If $u$ belongs to
$\D_{p,1}(\I_{\sigma+\e,p}^{1^-}),$ for some $\e>0,$ then $u$ is Stratonovitch integrable
on $[0,T]$ for any $T\in[0,1],$ and 
\begin{equation}
  \label{eq:14}
  \int_0^T u(s) \circ \d X_s =\div(\K^*_Tu)+\int_0^T \tilde{D}_Tu(s)\d
  s+ u(T)X(T).
\end{equation}
\end{thm}
\begin{proof}
 For the latest sum of $R_T^\pi(u-u(T)),$ we have according to Theorem~\ref{thm:existence_trace-}, 
\begin{multline*}
\esp{\left|   \sum_{t_i\in\pi} \frac{1}{\delta_i}\int_{t_i}^{t_{i+1}}\!\!\int_{t_i}^{t_{i+1}}
    \K^*_T\n_r(u-u(T))(s)\d s \d r- \int_0^T\tilde{D}_Tu(r)\d r\right|^p}\\
  \begin{aligned}
 \le &\,c\,   \esp{
    \sum_{t_i\in\pi} \frac{1}{\delta_i} \int_{t_i}^{t_{i+1}}\!\!\int_{t_i}^{t_{i+1}}
    |\K^*_T(\n_r(u-u(T))(s)-\tilde{D}_Tu(r)|^p\d s\d r}\\
\le & \,c\, \esp{
    \sum_{t_i\in\pi} \frac{1}{\delta_i}
    \int_{t_i}^{t_{i+1}}\!\!\int_{t_i}^{t_{i+1}} |s-r|^{p\eta}\lVert
    \n_r u\rVert^p_{\I^{1^-}_{\sigma+\e,p}}\d s\d r}\\
\le & \,c\,|\pi|^{p\eta} \lVert
    u\rVert^p_{\D_{p,1}(\I^{1^-}_{\sigma+\e,p})}.    
  \end{aligned}
\end{multline*}
Therefore, the latest sum of $R_T^\pi(u-u(T))$ converges in
$L^p(\Omega)$ (and thus in 
probability) to $\int_0^T \tilde{D}_Tu(s)\d s.$
In order to conclude, note that in virtue of the continuity of the
divergence, the first term of $R_T^\pi(u-u(T))$ tends to
$\div(\K^*_T(u-u(T))),$ see \cite{nualart.book}. 
\end{proof}
\begin{lemma}
  \label{lem:terme_additionnel}
Under the assumptions of Theorem \ref{thm:strato_int_inf}, for any 
$0\le S\le T\le 1,$ $u\car_{[0,S]}$ is Stratonovitch integrable on
$[0,T]$ and  we have
\begin{equation}
\label{eq:restriction_integrale}
\int_0^T (u(r) -u(S))\car_{[0,S]}(r) \circ \d X_r =\int_0^S u(r) \circ \d X_r, 
\end{equation}
for any $0\le S\le T\le 1.$
\end{lemma}
\begin{proof}
  According to Eqn. (\ref{eq:deltaaxi}) and to Lemma \ref{lem:triangularite}, we have
  \begin{equation*}
    R_T^\pi(p_S(u-u(S)))=R_S^\pi(u-u(S))+u(S)X^\pi(S).
  \end{equation*}
According to Theorem \ref{thm:strato_int_inf}, the right-hand-side sum
converges so that  $u\car_{[0,S]}$ is Stratonovitch integrable on
$[0,T]$ and Eqn. (\ref{eq:restriction_integrale}) follows by remarking
that $p_S(u-u(S))(T)=0.$
\end{proof}
\begin{rem}
  For the hypothesis `` $\K^*_1$ is continuous from
$\I_{\sigma,p}^{1^-}$ into $\Hol(\eta)$'' to hold, in view of the
examples cited above, this requires that $\sigma$ to be greater than $1/2-\alpha+1/p+\eta.$
\end{rem}
For $\alpha>1/2,$ the map $\K$ is still a regularizing operator so
that the hypothesis are much weaker. Following the very same lines, we
can prove:
\begin{thm}
\label{thm:existence_trace+}
Let $\alpha >1/2.$ Assume that Hypothesis \ref{A1} holds. Assume furthermore that $\K^*_1$ is
continuous from $\L^p$ into $\I^-_{\alpha-1/2,p}$ for some
$p>(\alpha-1/2)^{-1}.$ If $u$ belongs to $\D_{p,1}(\L^p),$ then, for
any $T\in[0,1],$  there exists a measurable and
integrable process, denoted by $\tilde{D}_Tu$ such that, for almost any~$r,$
\begin{equation*}
 \label{eq:15}
\esp{|\n_r \K^*_Tu(s)-\tilde{D}_Tu(r)|^p}^{1/p}\le \,c\,|s-r|^{\alpha-1/2-1/p} \lVert \n_r u\rVert_{L^p(\Omega\times[0,1])}.  
\end{equation*}
Moreover,
  \begin{equation*}
  \label{eq:16}
    \esp{\lVert \int_0^.\tilde{D}u(r)\d r\rVert_{\Hol(1-1/p)}^p}\le c \lVert
    u\rVert_{\D_{p,1}(\L^p)}^p.
  \end{equation*}
\end{thm}
\begin{thm}
\label{thm:strato_integ_sup}
Assume that Hypothesis \ref{A1} holds for $\alpha>1/2.$ Assume furthermore that $\K^*_1$ is
continuous from $\L^p$ into $\I^-_{\alpha-1/2,p}$ for some
$p>(\alpha-1/2)^{-1}.$ If $u$ belongs to $\D_{p,1}(\L^p),$ then  for
any $T\in[0,1],$ $u$ is Stratonovitch integrable on $[0,T]$ and 
\begin{equation*}
  \label{eq:def_intstrato}
  \int_0^T u_s \circ \d X_s=\div(\K^*_Tu)+\int_0^T \tilde{D}_Tu(s)\d s.
\end{equation*}
\end{thm}
\begin{rem}
The difference in this case is that $\L^p([0,1])$ is stable by the
maps $p_T$ so that we immediatly have:
  \begin{equation*}
    \label{eq:restriction_integrale2}
    \int_0^T u(s) \circ \d X_s=\int_0^1 u(s)\car_{[0,T]}(s)\circ \d X_s,
  \end{equation*}
in both theorems \ref{thm:strato_int_inf} and \ref{thm:strato_integ_sup}.
\end{rem}
Coming back to $\text{SS}_{\pi}(u),$ we have:
\begin{multline*}
\text{SS}_{\pi}(u)=\div\biggl(\sum_{t_i\in
  \pi}\frac{1}{\delta_i}\int_{t_i}^{t_{i+1}}\kern-8pt  u_s\, \d s\
\bigl(K(t_{i+1},.)-K(t_i,.)\bigr)\biggr)\\ +\sum_{t_i\in \pi} \frac{1}{\delta_i}
    \int_{t_i}^{t_{i+1}}
    \Bigl(K(\n.u_s)(t_{i+1})-K(\n.u_s)(t_{i})\Bigr)\, \d s
\end{multline*}
The trace-like term is
similar to those we had to treat in the previous theorems. The
difference is that its limit is formally $\int_0^1 (\K\n)_su(s)\d s$
instead of $\int_0^1 \n(\K^*_1u)(s)\d s$ in Theorems
\ref{thm:existence_trace-} and \ref{thm:existence_trace+}. We thus
need some regularity of the map $s\mapsto \n_s u(r)$ which is
something less easy to verify than properties on the map $s\mapsto
\n_r u(s).$ This restriction reduces the interest of this approach.
\begin{thm}
\label{thm:fin+}
Assume that Hypothesis \ref{A1} holds for $\alpha>1/2.$ Assume furthermore that $\K$ is
continuous from $\L^p([0,1])$ into $\I^-_{\alpha-1/2,p}$ for some
$p>(\alpha-1/2)^{-1}.$ If $u$ belongs to $\D_{p,1}(\L^p([0,1])),$ then there exists a measurable and
integrable process, denoted by $\hat{D}u$ such that, for almost any $r,$
\begin{equation}
\label{eq:def_de_D}
\esp{|(\K\n)_s u(r)-\hat{D}u(r)|}\le \,c\, |s-r|^\eta\lVert
D_. u(r)\rVert_{\D_{p,1}(\L^p([0,1]))}.
\end{equation}
Moreover, 
  \begin{equation}
 \label{eq:holder_de_D}
    \esp{\lVert \int_0^.\hat{D}u(r)\d r\rVert_{\Hol(1-1/p)}^p}\le c \lVert
    u\rVert_{\D_{p,1}(\L^p([0,1]))}^p.
  \end{equation}
Furthermore,  $\K^*_Tu$ belongs to $\Dom \div$ and the family
  $\text{SS}_\pi(u)$ converges in
  $L^2(\Omega)$ to $\div(\K^*_Tu) +\int_0^T \hat{D}u(s)\d s.$
\end{thm} 
\begin{rem}
  For $u$ belonging to $\D_{p,1}(\L^p([0,1]))$ and cylindric, it is
  easy to see that 
  \begin{equation}
    \label{eq:20}
    \int_0^1 \hat{D}u(r)\d r=\int_0^1 \hat{D}_1u(r)\d r.
  \end{equation}
According to (\ref{eq:holder_de_D}) and (\ref{eq:16}), this remains
true for any $u\in \D_{p,1}(\L^p([0,1])).$
\end{rem}
\begin{rem}
  For $\alpha<1/2,$ one could also state a similar theorem but it
  would  be practically of little use since it is rather hard to determine
  whether
  \begin{displaymath}
    \esp{\int_0^1 \lVert \n_.u(s)\rVert_{\S_{1+\eta-\alpha,2}}^2\d s}
    \text{ is finite.}
  \end{displaymath}
\end{rem}
\section{Regularity}
\label{sec:regul}
There are two kinds of regularity results which may be interesting :
continuity with respect to the time variable and continuity with
respect to the kernel. Actually, when one thinks to the generalized
fBm (see Example \ref{ex:gfbm}), the complete identification  of the model requires the perfect knowledge of the
function $H.$ Since that seems out of reach, one can naturally ask
how much an error on $H$ will modify the stochastic integral of a
given integrand.
The trace-like term  can be controlled via theorems
\ref{thm:existence_trace+} and \ref{thm:existence_trace-}. We are now
interested in the divergence part. 
 We denote by $\lVert
\K^*_{1}\rVert_{\alpha,p},$ the norm of $\K^*_1$ as a map from
$\I_{\alpha-1/2,p}^\prime$ into $\L^p.$
\begin{thm}
\label{thm:reg-}
Let  $\alpha\in(0,1/2)$ and $1< p <(1/2-\alpha)^{-1},$ assume that assumptions \ref{A1} and
\ref{A2} hold. Assume furthermore
that there exists $\e\in(0,1/p-(1/2-\alpha))$  such that $u$ belongs to
$\D_{p,1}(\I_{1/2-\alpha+\e,p}).$ Then, the process $\{\delta (\K^*_{t}u),\, t\in[0,1]\}$ admits a
modification with $\tilde\e$-H\"older continuous paths for any $\tilde{\e}<\e,$  and we have the
maximal inequality~:
\begin{equation*}
  \lVert \delta(\K^*_.u)\rVert_{L^p(\Omega; \I^+_{1/p^*+\tilde{\e},p^*})}\le \,c\, \lVert \K^*\rVert_{\alpha,p} \lVert
    u\rVert_{\D_{p,1}(\I_{1/2-\alpha+\e,p})}.
\end{equation*}
\end{thm}
\begin{proof}
Since  $1/2-\alpha+\e$ is strictly less
than $1/p,$  we know that for any
$T\in[0,1],$ $p_T u$ belongs to $\I_{1/p-\alpha+\e,p},$ see Proposition  \ref{prop:proprietes_int_RL}.
In view of Lemma \ref{lem:triangularite}, we have 
  $\delta(\K^*_tu)=\delta(\K^*_1(u\car_{[0,t]})).$ Therefore, for $g\in {\mathcal C}^\infty$ and $\psi$ a cylindric real-valued functional,
  \begin{multline*}
    \esp{\int_0^1 \delta\K^*_1(u\car_{[0,t]})g(t)\d t\, \psi}=
    \esp{\iint_{[0,1]^2} \K^*_1(u\car_{[0,t]})(r) g(t)\n_r\psi\d t\d r}\\
=\esp{\int_0^1 \K^*_1(uI^1_{1^-} g)(r)\n_r\psi \d r}
=\esp{\delta(\K^*_1(u.I^1_{1^-} g)\psi}.
  \end{multline*}
Thus,
\begin{equation}
  \label{eq:identification_deltakappastarug}
  \int_0^1 \delta(\K^*_tu)g(t)\d t=\delta(\K^*_1(u.I^1_{1^-}g))\text{ $\P$-a.s..}
\end{equation}
Since $p <(1/2-\alpha)^{-1},$ $1/2-\alpha<1/p,$  we can then apply Corollary \ref{cor:tambaca} with $t=\half-\alpha,$
  $r=1/p-\tilde{\e}$ and $s=\half-\alpha+\e.$
Since $g$ is deterministic,  we have
  \begin{equation}
    \label{eq:21}
    \lVert  \delta(\K^*_1(u.I^1_{1^-}g))\rVert_{L^p(\Omega)}\le c \lVert \K^*\rVert_{\alpha,p} \lVert
    u\rVert_{\D_{p,1}(\I_{1/2-\alpha+\e,p})}\lVert I^1_{1^-}g\rVert_{\I_{1/p-\tilde{\e},p}}.
  \end{equation}
We then obtain  that for $\psi\in
L^{p^*}(\Omega),$ for $g\in (\I^-_{1/p-1-\tilde{\e},p})^\prime,$
\begin{multline}
  \label{eq:6}
 |\esp{\int_0^1 \delta\K^*_1(u\car_{[0,t]})g(t)\d t\, \psi}|\\ \le  c \lVert
 \K^*\rVert_{\alpha,p} \lVert
 \psi\rVert_{L^{p^*}(\Omega)}\lVert g \rVert_{(\I^-_{1/p-1-\tilde{\e},p})^\prime}  \lVert
    u\rVert_{\D_{2,1}(\S_{1/2-\alpha,p})}.
\end{multline}
It follows that $\{\delta(\K^*_tu),\, t\in[0,1]\}$ belongs to
$(L^{p^*}(\Omega; \I^-_{-1+1/p-\tilde{\e},p}))^\prime,$ which is isomorphic to
$L^p(\Omega;\I^+_{1-1/p+\tilde{\e},p^*}),$ and that 
\begin{equation*}
  \lVert \delta(\K^*_.u)\rVert_{L^p(\Omega; \I^+_{1/p^*+\tilde{\e},p^*})}\le \,c\, \lVert \K^*\rVert_{\alpha,p} \lVert
    u\rVert_{\D_{p,1}(\I_{1/2-\alpha+\e,p})}.
\end{equation*}
This induces that there exists a modification of $\{\delta(\K^*_tu),\,
t\in[0,1]\}$ with $\tilde{\e}$-H\"older continuous sample-paths.
\end{proof}
\begin{rem}
  Note that $1$ belongs to $\I_{1/2-\e,2}$ for any $\e>0,$ thus we
  retrieve that $X_t=\div(\K^*_1p_tu)$ has a version with $(\alpha-\e)$-H\"older continuous sample-paths.
\end{rem}
If $\e>1/p-1/2+\alpha,$ we cannot apply  Lemma \ref{lem:tambaca} any more,
since $s=\half-\alpha+\e$ would be greater than $1/p.$ This is more
than a technical problem: in this situation, i.e., $u\in
\I_{\e+\half-\alpha,p},$ $u$ is continuous and $p_Tu$ does not
necesssary belongs to  $\I_{\e+\half-\alpha,p},$ so that the whole
principle of the above proof fails. However, as  Lemma
\ref{lem:inclusion_restriction} shows, if we consider $p_T(u-u(T))$ instead
of $P_Tu,$ this function belongs to $\I_{\e+\half-\alpha,p},$ for a
smaller $\e.$ Thus, we have:
\begin{thm}
\label{thm:regularite_divergence_-}
  Let  $\alpha\in (0,1/2)$ and $p> 1,$ assume that assumptions \ref{A1} and
\ref{A2} hold. Assume furthermore
that there exists $\e\in((1/p-1/2+\alpha)^+,1)$ such that $u$ belongs to
$\D_{p,1}(\I_{\e+\half-\alpha,p}^-).$ Then, for any $\tilde{\e}<\e,$ the process $\{\delta (\K^*_{t}(u-u(t))),\, t\in[0,1]\}$ admits a
modification with $\tilde{\e}$-H\"older continuous paths  and we have the
maximal inequality~:
\begin{equation}
\label{eq:11}
  \lVert \delta(\K^*_.(u-u(.)))\rVert_{L^p(\Omega; \I^+_{1/p^*+\tilde{\e},p^*})}\le \,c\, \lVert \K^*\rVert_{\alpha,p} \lVert
    u\rVert_{\D_{p,1}(\I_{\e+\half-\alpha,p}^-)}.
\end{equation}
\end{thm}
\begin{proof}
  Note that we are allowed to consider $u-u(t)$ since $1/p-1/2+\alpha<\e$
  implies that $\e+\half-\alpha>1/p$ and thus that
  $\I_{\e+\half-\alpha,p}^-$ is embedded in $\Hol(\e+\half-\alpha-1/p).$ The
  very same techniques as above show that 
  \begin{equation*}
    \int_0^1 \div(\K^*_t(u-u(t)))g(t)\d
    t=\div(\K^*_1(uI^1_{1^-}g-I^1_{1^-}(ug))),\text{ $\P$ a.s..}
  \end{equation*}
A classical integration by parts and then a fractional integration by
parts (see (\ref{int_parties_frac})) give that
\begin{equation*}
  \int_0^1 \div(\K^*_t(u-u(t)))g(t)\d
  t=-\div(\K^*_1(I^1_{1^-}(I_{0^+}^{-\zeta}u\, I_{1^-}^\zeta g))),\
  \P\text{ a.s..}
\end{equation*}
Now, we clearly have 
$$\lVert I^1_{1^-}( I_{0^+}^{-\zeta}uI_{1^-}^\zeta
g)\rVert_{\I_{1/2-\alpha,p}}=\lVert I_{0^+}^{-\zeta}u\, I_{1^-}^\zeta g\rVert_{\I_{-1/2-\alpha,p}}.$$
Applying Corollary  \ref{cor:tambaca} with $\zeta=1/2-\alpha+\e-1/p+\e^\prime,$ $t=-(1/2+\alpha),$
  $s+\zeta=1/2-\alpha+\e$ and
$r+s=t+1/p+\e^\prime$ for some $\e^\prime>0$ sufficiently small, we get
\begin{align*}
  \lVert \K^*_1(I^1_{1^-}(I_{0^+}^{-\zeta}uI_{1^-}^\zeta
  g))\rVert_{\L^p}& \le \,c\,   \lVert
  I_{0^+}^{-\zeta}u\rVert_{\I_{s,p}^-}\lVert I_{1^-}^\zeta g\rVert_{\I_{r,p}^-}  \\
&=\,c\,   \lVert
  u\rVert_{\I_{s+\zeta,p}^-}\lVert  g\rVert_{\I_{r-\zeta,p}^-}\\
&=
 \,c\, \lVert
  u\rVert_{\I_{1/2-\alpha+\e,p}^-}\lVert g\rVert_{\I_{-1+1/p-\e+\e^\prime,p}^-}.
\end{align*}
It follows as in the previous proof that $\{\div(\K^*_t(u-u(t))),\, t\ge
0\}$ belongs to $L^p(\Omega;\, \I_{1/p^*+\tilde{\e},p^*}^+)$ (with
$\tilde{\e}=\e-\e^\prime$) and that the
maximal inequality (\ref{eq:11}) holds.
\end{proof}
\begin{thm}
\label{thm:reg+}
For any $\alpha\in [1/2,1),$ assume that assumption \ref{A1} holds.
Let $u$ belong to $\D_{p,1}(\L^p)$ with $\alpha
p>1.$ The process $\{\delta (\K^*_{t}u),\, t\in[0,1]\}$ admits a
modification with $(\alpha-1/p)$-H\"older continuous paths and we have
the maximal inequality~:
\begin{equation*}
  \label{eq:2}
  \lVert \delta (\K^*_{.}u)\rVert_{L^p(\Omega;\Hol(\alpha-1/p))}\le c
  \lVert \K^*_{1}\rVert_{\alpha,2} \lVert u \rVert_{\D_{p,1}}.
\end{equation*}
\end{thm}
\begin{proof}
  We begin as in  Theorem \ref{thm:reg-}  until Eqn.
(\ref{eq:identification_deltakappastarug}). Since $\alpha>\half,$ it is clear that $\K$ is
  continuous from $\L^2([0,1])$ into $\I_{\alpha-\half,2}$ thus that
  $\K^*$ is continuous from $\I_{\alpha-\half,2}^\prime$ in
  $\L^2([0,1]).$ Since $\I_{\alpha-\half,2}$ is continuously embedded
  in $\L^{(1-\alpha)^{-1}},$ it follows that
  $\L^{1/\alpha}=(\L^{1/(1-\alpha)})^\prime$ is continuously embedded in
  $\I_{\half-\alpha,2}.$ Since $ u$ belongs to $\D_{p,1}(\L^p),$ the
  generalized H\"older inequality implies that
\begin{equation*}
\lVert   uI^1_{1^-}g  \rVert_{\L^{1/\alpha}}\le 
   \lVert   u \rVert_{\L^p} \lVert   I^1_{1^-}g \rVert_{\L^{(\alpha-1/p)^{-1}}}.
\end{equation*} 
It follows that $\{\delta (\K^*_{t}u),\, t\in[0,1]\}$ belongs to
$L^p(\Omega;\I^+_{1,(1-\alpha+1/p)^{-1}})$ with
\begin{displaymath}
  \lVert \delta (\K^*_{.}u)\rVert_{L^p(\Omega;\I^+_{1,(1-\alpha+1/p)^{-1}})}\le c
  \lVert \K^*_{1}\rVert_{\alpha,2} \lVert u \rVert_{\D_{p,1}}.
\end{displaymath}
The proof is completed remarking that
$1-1/(1-\alpha+1/p)^{-1}=\alpha-1/p$ so that
$\I^+_{1,(1-\alpha+1/p)^{-1}}$ is embedded in $\Hol(\alpha-1/p).$
\end{proof}
\begin{rem}
  These results extend similar results in \cite{alos00} in the sense
  that the assumptions on the kernel and on the integrand are here
  much weaker for the same conclusion. 
\end{rem}
\section{Itô Formula}
\label{sec:ito}
We are now interested in non-linear transformations of Itô-like
processes:
\begin{equation}
  \label{eq:3}
  Z(t)=z+\int_0^t u(s)\circ \d X_s,
\end{equation}
for a sufficiently regular $u.$ 
The  It\^o formula for fBm-like processes has already a long
history. There are two technical barriers: it is relatively easy to
prove It\^o formula for $\alpha>1/2,$ since we then have a process
more regular than the ordinary Brownian motion and all the limiting
procedures are straightforward (cf. \cite{dai96,decreusefond95_5,decreusefond96_2}) . Harder is the situation where $\alpha$
belongs to $(0,\, 1/2],$  Al\`os et al. \cite{alos01} obtained a
formula for the
fBm of Hurst index greater than $1/4.$
 By a very different procedure, Gradinaru et al. \cite{russo02} were
 able to include $1/4$ in the domain of validity of the formula. In
 another  different approach, Feyel et  al. \cite{MR1873303} also gave a formula for any
 Hurst index via analytic continuation of the formula obtained for
 $\alpha\ge 1/2.$  Carmona
 et al. \cite{coutin02} obtained an  It\^o formula for $\alpha >1/6,$
 for a class of processes similar to our so-called Volterra processes.
 
The following results owes much to the paper  \cite{coutin02} which
shows that it was possible to go beyond the barrier $1/4,$ to the
paper \cite{alos00} which gives the simplest expression of the It\^o
formula and to the work \cite{russo02} which emphasizes the importance
of symmetrization. Actually, the key remark is that there exists integrands $u$ for which 
\begin{equation}
\begin{split}
\label{eq:32}
R_h(u)& :=
 h^{-1} \int_0^1
\Bigl(\K^*_1p_{t+h}u(s)-\K^*_1p_tu(s)\Bigr)\Bigl(\K^*_1p_{t+h}u(s)+\K^*_1p_tu(s)\Bigr)\d
 s    \\
&= h^{-1} \int_0^1 \ks\pth u(s)\ks (p_t+p_{t+h})u(s)\d s\\
&= h^{-1} \int_0^1 (\ks p_{t+h}u(s)^2- \ks p_{t}u(s)^2)\d s,
\end{split}
  \end{equation}
has a finite limit. If $u\equiv 1,$ since
$\car_{[0,t)}=I_1^*(\varepsilon_t),$ it follows from the definition of
$\K$ that $\ks p_{t}\car =K(t,.)$ and thus $R_h(\car)=
h^{-1}(R(t+h,t+h)-R(t,t)),$ where $R$ is the covariance kernel of $X.$
For instance, if  $X$ is
the fBm with stationary increments, this expression is proportional to $h^{-1}((t+h)^{2\alpha}-t^{2\alpha}).$
The different barriers can be explained from the behavior of this last
term, whose limit is clearly $t^{2\alpha-1}.$ When
$\alpha>1/2,$ this is a bounded function of $t$ so easily controllable
in the limiting procedures. For $\alpha \in (1/4,1/2),$ it is no
longer bounded but still in $\L^2([0,1]).$ When,  $\alpha
<1/4,$  we only have an $\L^p$
integrable function for $1-p^{-1}<2\alpha.$
\begin{hyp}\label{A3}
  Let ${\mathcal R}$ the set of processes such that $R_h(u),$ as
  defined in (\ref{eq:32}), has a finite limit in $L^1(\Omega).$ We
  assume that $\ks$ is such that ${\mathcal R}$ is non-empty.
\end{hyp}
\begin{lemma}
  Let $\alpha\in(0,1),$ be given and assume that hypothesis \ref{A1},
  \ref{A2} and\ref{A3} hold. Let $u$ be a cylindric process, belonging
  to $\mathcal R.$ Let 
$$n_\alpha:=\inf\{n:\ 2n\alpha>1\}.$$
For any $f\in {\mathcal C}_b^{n_\alpha},$ i.e., $n_\alpha$-times
differentiable with bounded derivatives, we have
\begin{equation}
  \label{eq:33}
  \begin{split}
    \frac{d}{dt}\esp{\vphantom{\int_0^1}f(Z_t)\psi}&=\esp{\vphantom{\int_0^1}f^\prime(Z_t)(\K\n)_t(u(t)\psi)}\\
&+\frac{1}{2}\,\esp{ f^{\prime\prime}(Z_t)\psi \frac{d}{dt}\int_0^1
  \ks(p_t u)(s)^2\d s}\\
&+\esp{u(t)f^{\prime\prime}(Z_t)\psi (\K\n)_t\Bigl(\int_0^t
  (\K\n)_ru(r)\d r\Bigr)}\\
&+\esp{\vphantom{\int_0^1}u(t)f^{\prime\prime}(Z_t)\delta\bigl( (\K\n)_t(\ks p_t
  u)\bigr)\psi }.
  \end{split}
\end{equation}
\end{lemma}
\begin{proof}
  Introduce the function $g$ as
\begin{displaymath}
  g(x)=f(\frac{a+b}{2}+x)-f(\frac{a+b}{2}-x).
\end{displaymath}
This function is even, satisfies
$$g^{(2j+1)}(0)=2f^{(2j+1)}((a+b)/2) \text{ and
  }g(\frac{b-a}{2})=f(b)-f(a).$$
Applying the Taylor formula to $g$
between the points $0$ and $(b-a)/2,$ we get
\begin{multline*}
  f(b)-f(a)=\sum_{j=0}^{n-1}
  \frac{2^{-2j}}{(2j+1)!}(b-a)^{2j+1}f^{(2j+1)}(\frac{a+b}{2})\\+
  \frac{(b-a)^{2n}}{(2n)!}\int_0^1 \lambda^{2n-1}g^{(2n)}(\lambda a
  +(1-\lambda)b)\d \lambda.
\end{multline*}
We thus have
\begin{multline}
\label{eq:25}
\esp{(f(Z_{t+h})-f(Z_t))\psi}=\sum_{j=0}^{n_\alpha-1}\frac{2^{-2j}}{(2j+1)!}\esp{(b-a)^{2j+1}f^{(2j+1)}(\frac{a+b}{2})\ \psi}\\
+\frac{1}{2n_\alpha!}\esp{(Z_{t+h}-Z_t)^{(2n_\alpha)}\int_0^1
  r^{2n_\alpha-1}g^{(2n_\alpha)}(r Z_t +(1-r)Z_{t+h})\d r\ 
  \psi}.
\end{multline}
We need to prove that, when divided by $h,$ the latter quantity has a
limit when $h$ goes to $0.$ It turns out that the sole contributing
term is the first one.
We first show that $n_\alpha$ is chosen sufficiently large so that  the last term
vanish. Since $Z$ belongs
$L^2(\Omega;\, \Hol(\alpha-\e))$ for any $\e>0,$ and since $g^{(2n_\alpha)}$
is bounded, the last term is bounded by a constant times
$h^{2n_\alpha(\alpha-\e)}.$ Hence, this last term divided by $h$ vanishes when $h$
goes to $0.$
We next deal with the first order term.  Since $u$ is cylindric,
\begin{equation}
\label{eq:26}
  Z_t=\delta(\ks p_t u)+\int_0^t \ks(\n_sp_tu)(s)\d s.
\end{equation}
Substitute Eqn. (\ref{eq:26}) into the first order term and use
integration by parts formula, this yields to:
\begin{multline*}
  \esp{(Z_{t+h}-Z_t)f^\prime(\xth)\psi}\\
  \shoveleft{=\esp{\int \ks(\pth u)(s)
      \n_s(f^\prime(\xth)\psi)\d s}}\\
  \shoveright{+ \esp{f^\prime(\xth)\psi\int \ks(\pth\n_s u)(s)\d s\ }}\\
  \shoveleft{=\esp{f^\prime(\xth)\int  \ks(\pth u)(s)  \n_s\psi\d s\ }}\\
  \shoveright{+\esp{f^{\prime\prime}(\xth)\psi\int \ks(\pth u)(s) \n_s(\xth)\d s\ 
      }}\\
  + \esp{f^\prime(\xth)\psi\int_0^1 \ks(\pth\n_s u)(s)\d s\ 
    }=\sum_{i=1}^3A_i.
\end{multline*}
We can write $A_1$ as
\begin{equation*}
  A_1=\esp{\int_t^{t+h} u(s) (\K\n)_s\psi\d s\
    f^\prime(\xth)\psi},
\end{equation*}
by dominated convergence, it is then easily shown that
\begin{equation}
  \label{eq:27}
  \lim_{h\to 0} h^{-1}A_1 =\esp{u(t)f^\prime(Z_t)(\K\n)_t\psi }.
\end{equation}
By direct calculations, since $u$ is cylindric, we have
\begin{displaymath}
  \int_0^1 \ks(\pth\n_s u)(s)\d s=\int_t^{t+h} (\K\n)_su(s)\d s,\text{ thus,}
\end{displaymath}
\begin{equation}
  \label{eq:28}
  \lim_{h\to 0} h^{-1}A_3=\esp{f^\prime(Z_t)\psi(\K\n)_tu(t)}.
\end{equation}
Expanding $\n_s(Z_t+Z_{t+h}),$ we obtain
\begin{multline*}
  2 A_2=\esp{f^{\prime\prime}(\xth)\psi\int_0^1 \ks(\pth u)(s)\ks(p_t
    u+p_{t+h}u)(s)\d s\ }\\
  \shoveleft{+\esp{f^{\prime\prime}(\xth)\psi\int_0^1 \ks(\pth u)(s)\div\Bigl(\ks(p_t+p_{t+h})\n_su)\Bigr)\d s\ }}\\
\shoveleft{  +\text{E}\left[\int_0^1 \ks(\pth u)(s) \n_s\Bigl(\int_0^1(p_t+p_{t+h})
    (\K\n)_ru(r)\d r\d r\Bigr)\d s\right.}\\
  \times
  \left.f^{\prime\prime}(\xth)\psi\vphantom{\int_0^t}\right]
=\sum_{i=1}^3B_i.
\end{multline*}
According to Hypothesis \ref{A3},
\begin{equation}\label{eq:29}
  \lim_{h\to 0}h^{-1}B_1=\esp{\frac{\d}{\d t}\int_0^1 \ks(p_t u)(s)^2\d s f^{\prime\prime}(Z_t)\psi}.
\end{equation}
It is rather clear that
\begin{equation}\label{eq:30}
  \lim_{h\to 0}h^{-1}B_3=2\esp{u(t)(\K\n)_t(\int_0^t (\K\n)_ru(r)\d r)f^{\prime\prime}(Z_t)\psi}.
\end{equation}
To deal with $B_2,$ we need to apply once more the integration by
parts formula. This gives,
\begin{multline*}
  B_2=
  \text{E}\left[\int_0^1\int_0^1 \ks(\pth \n_r u)(s)\n_s\ks(p_t+p_{t+h})u(r)\d
      s\d r\right.\\
\shoveright{\left. \times f^{\prime\prime}(\xth)\psi\vphantom{\int_0^1}\right]}\\
 \shoveleft{ +\text{E}\left[\int_0^1 \ks(\pth u)(s)\int_0^1 \ks((p_t+p_{t+h})\n_su)(r)\right.}\\
\left.\times \n_r(f^{\prime\prime}(\xth)\psi)\d r \d s
\vphantom{\int_0^1}\right].
\end{multline*}
It follows from this expression that
\begin{align}
  \lim_{h\to 0}h^{-1}B_2& =2\, \esp{\int_0^1 (\K\n)_t(\ks p_t
    u)(r)\n_r(u(t)f^{\prime\prime}(Z_t)\psi)\d r} \notag\\
\label{eq:31}
&= 2\, \esp{u(t)f^{\prime\prime}(Z_t)\psi\, \delta\bigl( (\K\n)_t\ks p_t
  u\bigr)}.\vphantom{\int_0^1}
\end{align}
The remaining terms are of the form
 \begin{multline*}
   \esp{(Z_{t+h}-Z_t)^{2j+1}f^{(2j+1)}(\xth)\psi}\\
   =\esp{\int_0^1 \ks(\pth
     u)(s)\n_s\Bigl((Z_{t+h}-Z_t)^{2j}f^{(2j+1)}(\xth)\psi\Bigr)\d s}\\
   + \esp{(Z_{t+h}-Z_t)^{2j}f^{(2j+1)}(\xth)\psi\int_0^1\ks(\pth\n_s u)(s)\d s\ 
     }=C_1+C_2.
 \end{multline*}
 By dominated convergence, it is clear that $h^{-1}C_2$ vanishes as
 $h$ goes to $0.$ As to $C_1,$ it can be splitted into three parts
\begin{multline*}
  C_1=2j\, \text{E}\left[(Z_{t+h}-Z_t)^{2j-1}f^{(2j+1)}(\xth)\psi\right.\\
\shoveright{\left.\times \int_0^1 \ks(\pth u)(s) \n_s(Z_{t+h}-Z_t)\d
    s\right]}\\
\shoveleft{+\text{E}\left[\vphantom{\int_0^1}(Z_{t+h}-Z_t)^{2j}f^{(2j+2)}(\xth)\psi\right.}\\
\shoveright{\left. \int_0^1 \ks(\pth
    u)(s) \n_s(\xth)\d s\right]}\\
\shoveleft{+\text{E}\left[\vphantom{\int_0^1}(Z_{t+h}-Z_t)^{2j}f^{(2j+1)}(\xth)\right.}\\
\left.\times \int_0^1 \ks(\pth u)(s)\n_s \psi\d s\right]=\sum_{i=1}^3 D_i.
\end{multline*}
By dominated convergence, $h^{-1}D_3$ vanishes as $h$ goes to $0.$
Expanding the Gross-Sobolev derivative  $D_2,$ we get
\begin{multline*}
  2
  D_2=\text{E}\left[\vphantom{\int_0^1}f^{(2j+2)}(\xth)(Z_{t+h}-Z_t)^{2j}\psi \right.\\
\shoveright{\left.\times\int_0^1 \ks(\pth u)(s)\ks(p_t
    u+p_{t+h}u)(s)\d s \right]\hphantom{.}}\\
 \shoveleft{ +\text{E}\left[\vphantom{\int_0^1}f^{(2j+2)}(\xth)(Z_{t+h}-Z_t)^{2j}\psi\right.}\\
\shoveright{\left.\times\int_0^1 \ks(\pth u)(s)\div\Bigl(\ks(p_t\n_s
    u+p_{t+h}\n_su)\Bigr)\d s  \right]\hphantom{.}}\\
\shoveleft{   +\text{E}\left[f^{(2j+2)}(\xth)(Z_{t+h}-Z_t)^{2j}\psi\vphantom{\int_0^t}\right.}\\
  \times
  \left.\int_0^1 \ks(\pth u)(s) \n_s\Bigl(\int_0^1(p_t+p_{t+h})(\K\n)_ru(r)\d r\Bigr)\d s\right].
\end{multline*}
Following the reasoning applied to $A_2,$ we see that all the terms in
the integrals are converging a.s. (when divided by $h$) to a finite
limit, since there still is a factor $(Z_{t+h}-Z_t)^{2j},$ with $j>0,$
the product converges to $0.$ By dominated convergence, the
convergence can be seen to hold in $L^1(\Omega),$ thus $h^{-1}D_2$ goes
to $0$ as $h$ goes to $0.$ The really difficult term is $D_1.$ For the
sake of clarity, we only treat the case $j=1.$ For $j=1,$
\begin{multline*}
  D_1=\text{E}\left[\vphantom{\int_0^1}(Z_{t+h}-Z_t)f^{(3)}(\xth)\psi\right.\\
\shoveright{\left.\times\int_0^1 \ks(\pth u)(s) \n_s(Z_{t+h}-Z_t)\d
    s\right]\hphantom{.}}\\
\shoveleft{  =2\ \esp{(Z_{t+h}-Z_t)f^{(3)}(\xth)\psi\int_0^1 \ks(\pth u)(s)^2\d s}}\\
\shoveleft{  +2\ \text{E}\left[\vphantom{\int_0^1}(Z_{t+h}-Z_t)f^{(3)}(\xth)\psi\right.}\\
\shoveright{\left.\times\int_0^1 \ks(\pth u)(s) \div(\ks((p_t+p_{t+h})\n_su))\d
    s\right]\hphantom{.}}\\
\shoveleft{  +2\ \text{E}\left[\vphantom{\int_0^1}(Z_{t+h}-Z_t)f^{(3)}(\xth)\psi\right.}\\
\left.\times\int_0^1 \ks(\pth u)(s) \int_0^1
    \ks(\n^{(2)}_{r,s}(p_t+p_{t+h})u)(r)\d r\d
    s\right].
\end{multline*}
Dominated convergence implies that the last term, divided by $h,$
vanishes as $h$ goes to $0.$ For the two other summands, the idea is
always the same, each time there is a divergence term, we apply
integration by parts formula. Then, each new term is treated by the
previous methods. For instance, the most difficult term to handle is
one of the term which comes from derivative of the divergence in the first summand:
\begin{multline*}
  \esp{f^{(3)}(\xth)\psi\int_0^1 \n_r(\int_0^1 \ks(\pth u)(s)^2\d s)\ks(\pth u)(r)\d
    r}\\
\shoveleft{ = \text{E}\left[\vphantom{\int_0^1}f^{(3)}(\xth)\psi\right.}\\
\shoveright{\left.\times\int_0^1 \int_0^1 \ks(\pth u)(s)\ks(\pth \n
    _ru)(s)\ks(\pth
    u)(r)\d r\d s\right]\hphantom{.}}\\
\shoveleft{ = \text{E}\left[\vphantom{\int_0^1}f^{(3)}(\xth)\psi\right.}\\
\left.\times \int_t^{t+h}u(s)\K\Bigl(\int_0^1\ks(\pth\n_ru)(.)\ks(\pth
    u)(r)\d r\Bigr)(s)\d s\right].
\end{multline*}
Once again, in this form, it is clear that this term, divided by $h,$
converges to $0.$ All the remaining term are treated likewise and do
not contribute. Thus from Eqn. (\ref{eq:33}) follows from 
(\ref{eq:27}), (\ref{eq:28}), (\ref{eq:29}), (\ref{eq:30}) and (\ref{eq:31}).
\end{proof}
Since $u$ is cylindric, all the terms of  (\ref{eq:33}) are integrable
with respect to $t,$ we thus have
\begin{cor}
  Under the assuptions of the previous lemma,
we have, 
\begin{align*}
\esp{  f(Z_t)\psi}& =\esp{f(x)\psi}+\esp{\int_0^t f^\prime(Z_s)(\K\n)_s(u(s)\psi)\d s}\\
&+\frac{1}{2}\esp{\psi\int_0^t f^{\prime\prime}(Z_s)\frac{d}{ds}\int_0^1 \ks p_s
u(r)^2\d r\d s}\\
&+\esp{\psi\int_0^t u(s)f^{\prime\prime}(Z_s)(\K\n)_s\Bigl(\int_0^s
(\K\n)_ru(r)\d r\Bigr)\d s}\\
&+ \esp{\psi\int_0^t u(s)f^{\prime\prime}(Z_s)\delta\Bigl(\ks p_s (\K\n)_su\Bigr)\d s},
\end{align*}
for any $\psi$ such that $\n\psi$ belongs to $\Dom\K.$
\end{cor}

Since $(\K\n)$ is a derivation operator, we obtain after a few
manipulations:
Since $\K\n$ is a derivation operator, we have
\begin{align*}
\esp{  f(Z_t)\psi}& =\esp{f(x)\psi}+\esp{\int_0^t (\K\n)_s(f^\prime(Z_s)u(s)\psi)\d s}\\
&+\frac{1}{2}\esp{\psi\int_0^t f^{\prime\prime}(Z_s)\frac{d}{ds}\int_0^1 \ks p_s
u(r)^2\d r\d s}\\
&-\esp{\psi\int_0^t u(s)f^{\prime\prime}(Z_s)\K\ks (p_s u)(s)\d s}.
\end{align*}
This means that for any $t,$ we have a.e.,
\begin{equation}\label{eq:36}
  \begin{split}
  f(Z_t)& =f(x)+\int_0^t f^\prime(Z_s)u(s)\circ \d X_s\\
&+\frac{1}{2}\int_0^t f^{\prime\prime}(Z_s)\frac{d}{ds}\int_0^1 \ks p_s
u(r)^2\d r\d s\\
&- \int_0^t u(s)f^{\prime\prime}(Z_s)\K\ks (p_s u)(s)\d s.   
  \end{split}
 \end{equation}
 \begin{rem}
   It has to be noted that in \cite{decreusefond02}, we announced an
   It\^o formula for general $u$ and any $\alpha\in (0,1).$ This is
   unfortunately wrong for $\alpha\in(0,1/2).$
 Actually, starting from (\ref{eq:36}), the problem is now to pass to the limit. For the very first term of
the righthandside of (\ref{eq:36}), we need to find a class of
processes $u$ for which $f\circ Z. u$ is Stratonovich integrable. The
most restrictive part is to find conditions under which this process
has a ``trace'' in the sense of Theorem \ref{thm:existence_trace-}. It
is important to note that 
\begin{multline*}
  \n_r Z_t=\ks p_t(u-u(t))(r)+\delta(\ks p_t \n_r (u-u(t)))+\n_r \int_0^t
  (\K\n)_s(u-u(t))(s)\d s\\
+X(t)\n_r u(t)+u(t)K(t,r) 
\end{multline*}
and thus,  we have
\begin{multline*}
    \K(\n_. Z_t)(r)=\K(\ks p_t(u-u(t)))(r)+\K(\delta(\ks p_t \n_. (u-u(t))))(r)\\+\K(\n_. \int_0^t
  (\K\n)_s(u-u(t))(s)\d s)(r)
+\K(X(t)\n_. u(t))(r)\\ +\K(u(t)K(t,.))(t). 
\end{multline*}
It is possible to impose hypothesis on $u$ such that the first four
terms of the previous equations have a signification when $r=t.$
Unfortunately, for the very last term, we have
\begin{equation*}
  \label{eq:37}
  \K(u(t)K(t,.))(t)=u(t)\frac{\partial }{\partial s}R(t,s)|_{s=t}.
\end{equation*}
In the case of the fBm with stationary increments, this is equal, up
to a constant, to $u(t) (s^{2\alpha-1}-(t-s)^{2\alpha-1})_{s=t}.$
Since this quantity is infinite for $\alpha<1/2,$  we haven't been able to go below $1/2.$ 
 \end{rem}
\begin{rem}If we don't have a trace term we can
  state the following result.
  \begin{thm}
     Let $\alpha\in(0,1),$ be given and assume that hypothesis \ref{A1},
  \ref{A2} and\ref{A3} hold. Let $u$ be a cylindric process, belonging
  to $\mathcal R.$ Let 
$$n_\alpha:=\inf\{n:\ 2n\alpha>1\}.$$ Let 
\begin{equation*}
  Z_t=\delta(\ks p_t u).
\end{equation*}
For any $f\in {\mathcal C}_b^{n_\alpha},$ i.e., $n_\alpha$-times
differentiable with bounded derivatives, we have
\begin{align*}
  f(Z_t)& =f(x)+\delta\Bigl(\K^*_t (u. f^\prime\circ Z)\, \Bigr)\\
&+\frac{1}{2}\int_0^t f^{\prime\prime}(Z_s)\frac{d}{ds}\int_0^1 \ks p_s
u(r)^2\d r\d s\\
&+ \int_0^t u(s)f^{\prime\prime}(Z_s)\delta\Bigl(\ks p_s (\K\n)_su\Bigr)\d s,
\end{align*}
for any $t,$ a.s..
  \end{thm}
  \begin{proof}
    The proof is exactly the same as the previous one.
  \end{proof}
If $u\equiv 1,$ we get the same result as in
\cite{alos00,coutin02,decreusefond96_2,MR1873303} valid for any $\alpha\in (0,1).$ If
$\K=\text{Id},$ i.e., $X$ is an ordinary Brownian motion, and $u$ is
not necessarily adapted, this formula coincides with that given  in \cite{ustunel88_1}.
\end{rem}
\section{Skorohod integral}
Since the term $\int_0^T \tilde{D}_Tu(s)\d s$ is a trace-like term, it
is reasonable to introduce the following definitions.
We now introduce a stochastic integral defined 
\begin{defn}
  We denote by $\Dom \div_{\K^*}, $ the set of processes $u$ belonging
  a.s. to $\Dom \K^*$ and such that $\K^*u$ belongs to $\Dom \div.$ We
  denote by $\Dom \div_X, $ the set of processes $u$ in $\Dom
  \delta_{\K^*}$ such that $\n \K^*u$ is $\P$-a.s. a trace class
  operator.
\end{defn}
\begin{defn}
  For $u\in \Dom\div_X,$ we define the stochastic integral of $u$ with
  respect to $X$ by
  \begin{align*}
    \int_0^1 u_s * \d X_s & \egaldef \int_0^1 (\K^*u)(s) \delta B_s +
    \trace(\n(\K^*u))
  \end{align*}
\end{defn}
To define the integral of $u$ between time $0$ and $t,$ we use Lemma
\ref{lem:triangularite}:
\begin{defn}
  For $u\in \Dom\div_X,$ we define the stochastic integral of $u$ with
  respect to $X$ between $0$ and $t$ by
  \begin{align*}
    \label{eq:definition_integrale0t}
    \int_0^t u_s * \d X_s & = \int_0^1 (p_tu)(s) *\d X_s\nonumber\\
    & = \int_0^t (\K^*_t u)(s) \delta B_s + \trace(p_t\n(\K^*_tu)),
\end{align*}
where the second equality follows by (\ref{eq:5}).
\end{defn}
Eqn. (\ref{eq:20}) has its equivalent in this setting :
\begin{lemma}
\label{lem:egalitetraces}
Assume that \ref{A1} and\ref{A2} hold. Let $u\in \Dom \K^*$ belong to $\D_{2,1}(\L^2([0,1]))$ and be such that $\n
u$ belong (a.s.) to $\Dom \K.$ Then $\trace(\n(\K^*u))$ is finite iff
$\trace((\K\n)u)$ is finite and they are equal.
\end{lemma}
\begin{proof}
  Since $\Dom \K^*\cap \Dom\K$ is a dense subset of $\L^2$, one can
  find $\{h_i, \, i\ge 1\}$ an ONB of $\L^2$ where for any $i,$ $h_i$
  belongs to $\Dom \K^*\cap \Dom\K.$ Set $\pi_n$ the orthogonal
  projection in $\L^2$ onto the vector space spanned by $h_1,\ldots, h_n.$ Let $V_k=\sigma\{\delta h_i,
  i=1,\ldots,k\}$ and consider $u_{k,n}=\pi_n\esp{P_{1/k}u \/ V_k}$
  where $P_t$ denote the Ornstein-Uhlenbeck semi-group of the Wiener
  process $X.$ It is known, see \cite[Lemma B.6.1]{ustunel2000}, that
  $u_k$ can be written as
\begin{equation*}
  u_{k,n}=\sum_{i=1}^nf_i^n(\delta h_1,\ldots, \delta h_k)h_i \text{ where } f_i\in {\mathcal
    C}^\infty \text{ for any } i,
\end{equation*}
and that $u_{k,n}$ converges to $u$ in $\D_{2,1}.$ Furthermore, it is
clear that we have
\begin{align}
  \trace((\K \n) u_{k,n})&=\trace\sum_{i,j}\partial_jf_i^n(\delta
  h_1,\ldots, \delta h_k)h_i\otimes \K h_j\nonumber\\
  &=\sum_{i,j}\partial_jf_i^n(\delta
  h_1,\ldots, \delta h_k) \int_0^1 h_i(s) (\K h_j)(s)\, ds\nonumber\\
  &=\sum_{i,j}\partial_jf_i^n(\delta
  h_1,\ldots, \delta h_k) \int_0^1 (\K^*h_i)(s)  h_j(s)\, ds\nonumber\\
  &=\trace(\n (\K^*u_{k,n})).\label{eq:egalitetraces}
\end{align}
Moreover, if $\trace((\K \n)u)$ exists a.s., then the series
\begin{equation*}
  \sum_i <(\K\n)u,\, h_i\otimes h_i>_{\L^2\otimes \L^2} \text{ is convergent.}
\end{equation*}
Thus, by Cauchy-Schwarz inequality,
\begin{multline*}
  \Bigl|\trace((\K \n) u_{k,n}) - \trace((\K \n)u)\Bigr| \\
  \le \sum_{i\le n} <(\K \n) u_{k,n}-(\K \n) u, h_i\otimes
  h_i>_{\L^2\otimes \L^2}+\sum_{i>n}|<(\K\n)u,\, h_i\otimes h_i>_{\L^2\otimes \L^2}|\\
  \le n. \lVert (\K \n)(u-u_{k,n})\rVert_{\L^2\ox \L^2} +
  \sum_{i>n}|<(\K\n)u,\, h_i\otimes h_i>_{\L^2\otimes \L^2}|.
\end{multline*}
As $n$ goes to infinity, the rightmost term converges a.s. to $0,$
hence for $\ee>0,$ one can find $n$ such that
\begin{equation*}
  \P(\sum_{i>n}|<(\K\n)u,\, h_i\otimes h_i>_{\L^2\otimes \L^2}>\ee/2)\le \ee/2.
\end{equation*}
Since $\K$ is a closed map, for this value of $n,$ one can find $k_n$
such that
\begin{equation*}
\P(  \lVert (\K \n)(u-u_{k_n,n})\rVert_{\L^2\ox \L^2}>\ee/2n)\le \ee/2.
\end{equation*}
For such $n$ and $k_n, $ we have
\begin{equation*}
  \P(\Bigl|\trace((\K \n) u_{k_n,n}) - \trace((\K \n)u)\Bigr| >\ee)\le \ee.
\end{equation*}
Hence there exists a subsequence $(k_j,n_j)$ such that $\trace((\K \n)
u_{k_j,n_j})$ converges $\P$-almost surely, thus that
$\trace(\n(\K^*u))$ is finite and that the two expressions are equal.
$\trace(\n(\K^*u))=\trace((\K \n)u).$

The very same reasoning holds true when $\trace(\n(\K^*u))$ is finite.
\end{proof}
Following \cite{nualart.book}, we know that when $u$ belongs to the
domain of the two integrals (that of definition \ref{def:strat-integr}
and that of the last definition), these two integrals  coincide.

A nice feature of this version of the stochastic integral is that we
can compute its transformation under absolutely continuous change of probability.
\begin{thm}
\label{cor:chgtproba}
Let $T(\omega)=\omega+Kv(\omega)$ be such that $v$ belongs to
${\mathbb D}_{p,1}(\L^2)$ for some $p>1$ and $T^*P \ll P$. Let $u$ be
such that $u$ and $u\circ T$ belong to $\Dom \delta_{\K^*}$ and $\n
\K^* u$ and $\n (\K^*u\circ T)$ are a.s. trace class operators. Then,
  \begin{displaymath}
\bigl(\int_0^1 u(s) *\d X_s\bigr)\circ T=   \int_0^1 (u\circ T)(s) * \d X_s + \int_0^1
\K^*(u\circ T)(s)  v(s)\, ds.
  \end{displaymath}
\end{thm}
\begin{proof}
  Theorem B.6.12 of \cite{ustunel2000} stands that
  \begin{equation*}
    \delta(\K^*u)\circ T= \delta(\K^*(u\circ T))+ \int
\K^*(u\circ T)(s)  v(s)\, ds+\trace((\n \K^*u)\circ T . \n v).
  \end{equation*}
  Proposition B.6.8 of \cite{ustunel2000} implies that
  \begin{equation*}
  \trace((\n \K^*u)\circ T . \n v)=\trace(\n (\K^*u\circ T))-\trace(\n
  \K^*u)\circ T.  
  \end{equation*}
  The proof is completed by substituting the latter equation into the
  former.
 \end{proof}
 For $u$ deterministic and $v$ adapted, this means that the law of the
 process $\{\int_0^t u_s dX_s - \int_0^t \K^*u(s) v(s)\, ds, t\ge
 0\},$ under $T^*P,$ is identical to the $P$-law of the process
 $\{\int_0^t u_s dX_s, t\ge 0\}.$

  \def\polhk#1{\setbox0=\hbox{#1}{\ooalign{\hidewidth\lower1.5ex\hbox{`}\hidewidth\crcr\unhbox0}}}
\providecommand{\bysame}{\leavevmode\hbox to3em{\hrulefill}\thinspace}
\providecommand{\MR}{\relax\ifhmode\unskip\space\fi MR }
\providecommand{\MRhref}[2]{%
  \href{http://www.ams.org/mathscinet-getitem?mr=#1}{#2}
}
\providecommand{\href}[2]{#2}

\end{document}